\documentclass[10pt, oneside, reqno]{amsproc}

\usepackage{amsmath, amssymb, amsfonts, amsthm, bm, graphicx}

\DeclareMathOperator{\conv}{conv}

\DeclareMathOperator{\supp}{supp}
\DeclareMathOperator{\im}{im}
\newcommand{\rma}{\boxplus}
\newcommand{\rmul}{\boxdot}
\newcommand{\subdiff}{\partial}

\newcommand{\fhi}{\varphi}
\newcommand{\Tr}{\mathsf{T}}
\newcommand{\norm}[1]{\lVert#1\rVert}

\newcommand{\M}{\mathbb{M}}
\newcommand{\R}{\mathbb{R}}
\newcommand{\one}{\mathbf{1}}
\newcommand{\vect}[1]{\bm{#1}}
\newcommand{\va}{\vect{a}}
\newcommand{\vb}{\vect{b}}

\newcommand{\ve}{\vect{e}}

\newcommand{\vj}{\vect{j}}
\newcommand{\vx}{\vect{x}}
\newcommand{\vy}{\vect{y}}

\newcommand{\vo}{\vect{o}}
\newcommand{\vp}{\vect{p}}
\newcommand{\vq}{\vect{q}}

\newcommand{\vv}{\vect{v}}

\newcommand{\collection}[1]{{\mathcal#1}}

\newcommand{\CP}{\collection{P}}

\newcommand{\CZ}{\collection{Z}}
\newcommand{\abs}[1]{|#1|}

\newcommand{\card}[1]{\lvert#1\rvert}

\theoremstyle{plain}
\newtheorem{theorem}{Theorem}
\newtheorem*{separationtheorem}{Separation Theorem}
\newtheorem{lemma}[theorem]{Lemma}
\newtheorem{corollary}[theorem]{Corollary}
\newtheorem{claim}{Claim}
\newtheorem{conjecture}[theorem]{Conjecture}
\newtheorem{proposition}[theorem]{Proposition}

\begin{document}

\bibliographystyle{amsplain}

\title[The local Steiner problem in normed spaces]{The local Steiner problem 
in finite-dimensional normed spaces}
\author{Konrad J. Swanepoel}
\thanks{This material is based upon work supported by the South African 
National Research Foundation under Grant number 2053752.}
\address{Department of Mathematical Sciences,
        University of South Africa, PO Box 392,
        Pretoria 0003, South Africa}
\email{\texttt{swanekj@unisa.ac.za}}
\date{\today}

\begin{abstract}
We develop a general method for proving that certain star configurations in 
finite-dimensional normed spaces are Steiner minimal trees.
This method generalises the results of Lawlor and Morgan (1994) that could 
only be applied to differentiable norms.
The generalisation uses the subdifferential calculus from convex analysis.
We apply this method to two special norms.
The first norm, occurring in the work of Cieslik, has unit ball the polar 
of the difference body of the $n$-simplex (in dimension $3$ this is the 
rhombic dodecahedron).
We determine the maximum degree of a given point in a Steiner minimal tree 
in this norm.
The proof makes essential use of extremal finite set theory.
The second norm, occurring in the work of Mark Conger (1989), is the sum of 
the $\ell_1$-norm and a small multiple of the $\ell_2$ norm.
For the second norm we determine the maximum degree of a Steiner point.
\end{abstract}

\maketitle

\section{Introduction}
There is a vast literature on Steiner minimal trees, mostly in graphs, Hamming 
space and related word spaces, the Euclidean plane, and the Manhattan plane, 
with applications in VLSI design \cite{KPS} and phylogenetics \cite{Cieslikphylo}.
See e.g.\ the monographs \cite{Cieslik2, HRW, IT}.
Cockayne \cite{Cockayne} considered Minkowski planes (two-dimensional 
normed spaces).
Cieslik initiated the study of Steiner minimal trees in general 
finite-dimensional normed spaces \cite{Cieslik, Cieslik4}.
Recently, other norms besides Euclidean and Manhattan have found applications 
\cite{CCKMWY, LS, LSD, LCLW, LX, LXZ}, while the geometry of certain 
high-dimensional normed spaces is related to the metric spaces found in 
the mathematical study of phylogenetic trees \cite{Cieslikphylo}.
These trees are also of interest in differential geometry; see especially 
the work of Morgan and his students \cite{Morgan, Morganbook, Alfetal, 
Conger, LM}.

\smallskip
We denote an $n$-dimensional normed space or \emph{Minkowski space} 
$(\R^n,\norm{\cdot})$ by $\M^n$ (see Section~\ref{mspaces}).
A \emph{Steiner tree} $T$ of a finite set of points $N$ in $\M^n$ is a 
tree with vertex set $V=V(T)$ and edge set $E=E(T)$ where $N\subseteq V$.
The points in $N$ are called \emph{nodes} or \emph{given points}, and the 
points in $V\setminus N$ \emph{Steiner points} or \emph{auxiliary points}.
The \emph{length} of a Steiner tree $T$ is
\[ \ell(T) := \sum_{\vx\vy\in E(T)}\norm{\vx-\vy}.\]
Any finite set of points $N$ in a Minkowski space has at least one tree 
of minimum length \cite{Cockayne}.
Such a shortest Steiner tree of $N$ is called a \emph{Steiner minimal tree} 
(SMT) of $N$.

A \emph{star} is an SMT such that one of the vertices, called the 
\emph{centre} of the star, (either a Steiner point or a node) is joined to 
all the other vertices.
The \emph{vertex figure} of a vertex in an SMT is the set of vectors from 
the centre of the star to the other vertices.
Note that only the directions of the vectors matter:
If a set of non-zero vectors is a vertex figure in some SMT, then we may 
replace each vector with another in the same direction, to form a new SMT.
In an earlier paper \cite{SwaN} we called the problem of describing the 
vertex figures of nodes and Steiner points the \emph{local Steiner problem}.
Here we give a solution of the local Steiner problem in Minkowski spaces, 
formulated in terms of norming functionals (Section~\ref{characterization}).

We denote the maximum degree of a given point in an SMT in $\M^n$, with 
the maximum taken over all possible SMTs, by $d(\M^n)$.
Let $s(\M^n)$ be the maximum degree of a Steiner point in an SMT in $\M^n$, 
where the maximum is taken over all possible SMTs in $\M^n$.
If $T$ is an SMT of $N$, then $T$ is clearly still an SMT of any $N'$ 
where $N\subseteq N'\subseteq V(T)$.
It follows that $s(\M^n)\leq d(\M^n)$.
We use the characterization of vertex figures to obtain information on 
these values for certain spaces.

Cieslik \cite{Cieslik}  proved that $d(\M^n)$ is bounded above by the 
\emph{Hadwiger number} or \emph{translative kissing number} of the unit 
ball $B(\M^n)$, i.e., the maximum number of mutually non-overlapping 
translates of $B(\M^n)$ that all touch $B(\M^n)$.
This already gives an upper bound depending only on $n$, since the Hadwiger 
number is bounded above by $3^n-1$ \cite{Hadwiger}.
He proved this more generally for all minimal spanning trees, a result 
rediscovered in \cite{RobSa} for the special case of $\ell_p$ norms.
Cieslik \cite{Cieslik}, \cite[Conjecture~4.3.6]{Cieslik2} also made the 
following conjecture:
\begin{conjecture}[Cieslik \cite{Cieslik, Cieslik2}]
The maximum degree of a given point in an SMT in any $n$-dimensional 
Minkowski space $\M^n$ satisfies \[d(\M^n)\leq 2^{n+1}-2,\] with equality 
if and only if $\M^n$ is isometric to the space $\M_Z^n$ which has unit 
ball $\conv([0,1]^n\cup[-1,0]^n)$.
\end{conjecture}
The unit ball of $\M_Z^2$ is an affine regular hexagon, and of $\M_Z^3$ 
an affine rhombic dodecahedron.
\begin{figure}
\begin{center}
\includegraphics[scale=0.42]{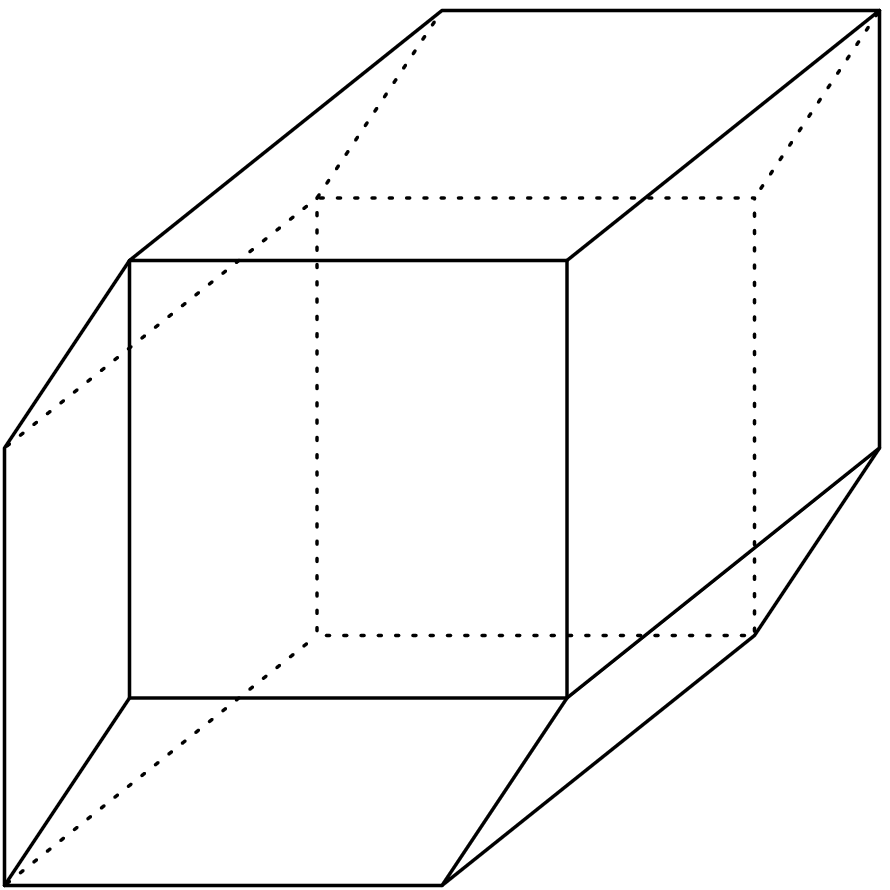}\qquad\qquad
\includegraphics[scale=0.5]{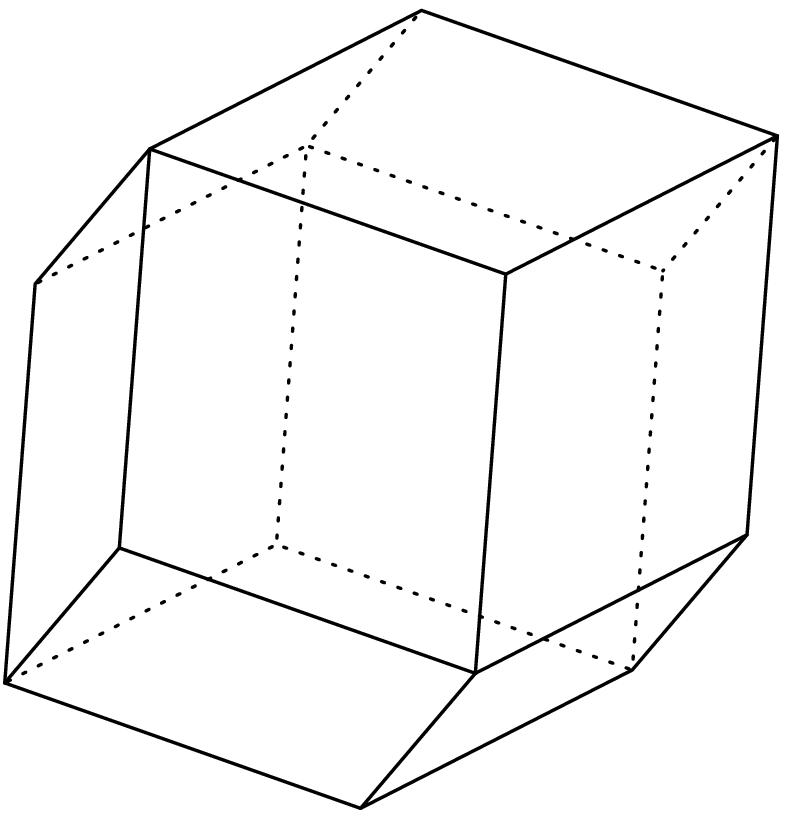}
\end{center}
\caption{Two representations of the rhombic dodecahedron}\label{fig3}
\end{figure}
The norm of $\M_Z^n$ is
\[ \norm{(\vx(1),\dots,\vx(n))} = \max_{\vx(i)\geq 0} \vx(i) - 
\min_{\vx(i)\leq 0} \vx(i) .\]
See Fig.~\ref{fig3} for two affine representations, one as in the conjecture, 
and the other with all faces congruent.
In Section~\ref{zonotopes} we represent $\M_Z^n$ in a different, 
more convenient way.

\begin{figure}[b]
\begin{center}
\includegraphics[scale=0.45]{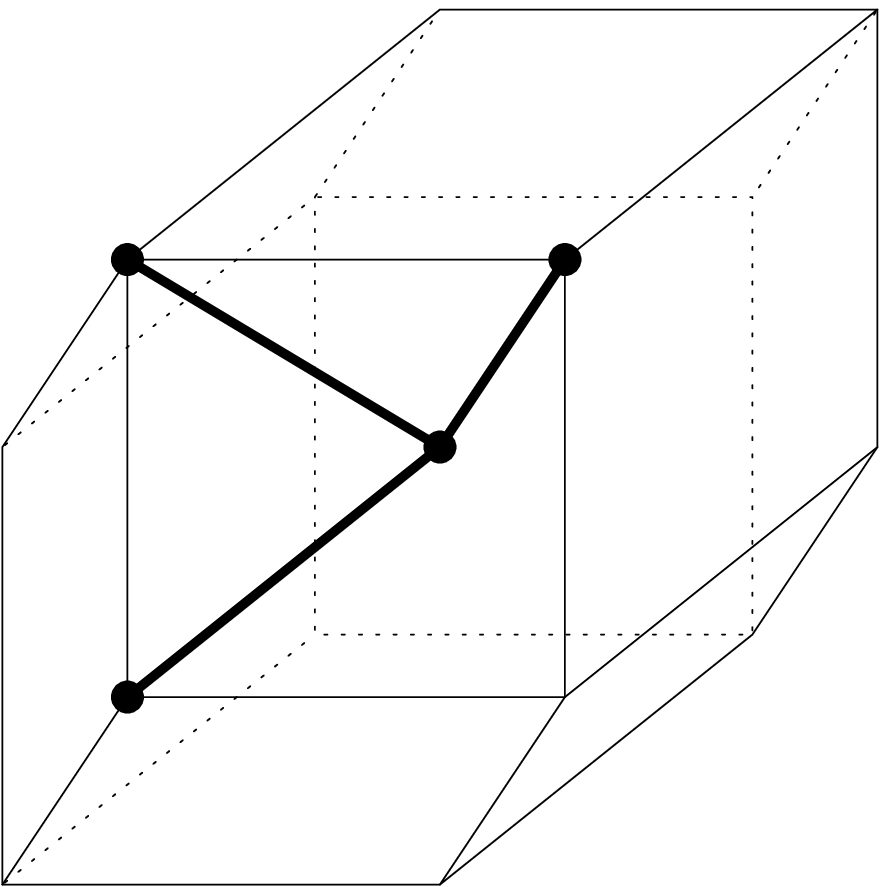}\qquad\qquad
\includegraphics[scale=0.45]{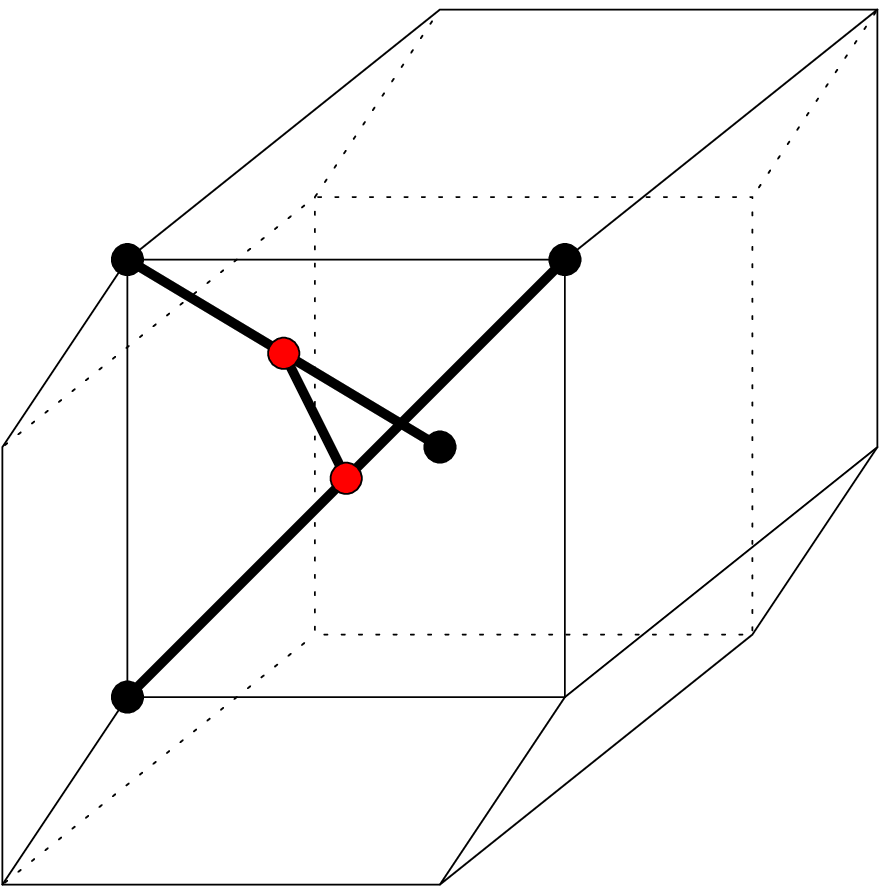}
\end{center}
\caption{Shortening a tree in the rhombic dodecahedral norm}\label{fig2}
\end{figure}

\begin{figure}
\begin{center}
\includegraphics[scale=0.53]{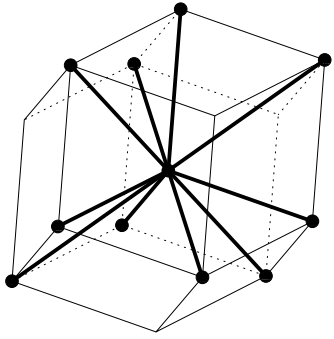}
\end{center}
\caption{A given point of degree $10$ in the rhombic dodecahedral 
norm}\label{fig4}
\end{figure}

Since any edge joining two points in an SMT can be replaced by a 
piece-wise linear path consisting of segments parallel to the vectors 
pointing to the extreme points of the unit ball, we obtain the following 
well-known lemma.
\begin{lemma}\label{polyhedral}
If the unit ball of $\M^n$ is a polytope with $v$ vertices, then 
$d(\M^n)\leq v$.
\end{lemma}
It follows that $d(\M_Z^n)\leq 2^{n+1}-2$.
Cieslik \cite{Cieslik4} proved the case $n=2$ of his conjecture.
However, we show that $d(\M_Z^n)<2^{n+1}-2$ for all $n\geq 3$, thus 
partially disproving the conjecture.
It is not difficult to see that $d(\M^3_Z)<14$ by shortening the star 
connecting the origin to the vertices of the rhombic dodecahedron.
If there exists an SMT with a given point of degree $14$, then the star 
joining $\vo$ to the vertices of the unit ball will be an SMT of the 
$14$ vertices together with the origin.
In Fig.~\ref{fig2} on the left we have three of the fourteen edges from 
the origin to the vertices.
These three vertices, one of degree $4$ and two of degree $3$, are on 
the same facet of the unit ball.
The length of this subtree is $3$.
On the right we show that it can be shortened to a length of $5/2$ by 
introducing two Steiner points.
This sub-configuration is therefore forbidden in any star joining the 
origin to a subset of the vertices, and it easily follows that 
$d(\M_Z^3)\leq 10$.
This turns out to be the correct value (Fig.~\ref{fig4}).
\begin{theorem}\label{thm1}
For any $n\geq3$, $d(\M^n_Z)=\binom{n+2}{\lfloor (n+2)/2\rfloor}$.
\end{theorem}
The proof is in Section~\ref{extremal}.
We also describe all possible vertex figures of a node in an SMT in 
$\M^n_Z$.
From the proof we obtain the following corollary.
\begin{corollary}\label{cor1}
In the Minkowski space $\M_Z^n$, the star joining $\vo$ to any set 
$A$ of non-zero points is an SMT of $A\cup\{\vo\}$ if and only the 
star joining $\vo$ to any subset $B$ of $A$ of size at most $4$ is 
an SMT of $B\cup\{\vo\}$.
\end{corollary}
For all $n\geq 2$ there is, up to isometries, a unique configuration 
attaining the maximum degree (Theorem~\ref{thm11} in Section~\ref{extremal}).
These spaces give the largest known degrees of SMTs in Minkowski spaces 
of dimensions $2$ to $6$, in fact $d(\M^3_Z)=10$, $d(\M^4_Z)=20$, 
$d(\M^5_Z)=35$, and $d(\M^6_Z)=70$.
For $n\geq 7$, $d(\ell_\infty^n)=2^n$ is larger.
It is not at all clear whether $\M^n_Z$ maximises $d(\M^n)$.

\bigskip
Morgan \cite[p.~42]{Morgan}, \cite[p.~129]{Morganbook} asks the 
following question:
\begin{conjecture}[Morgan \cite{Morgan, Morganbook}]
The maximum degree of a Steiner point in an SMT in any $n$-dimensional 
Minkowski space $\M^n$ satisfies \[s(\M^n)\leq 2^n.\]
\end{conjecture}
The asymptotically best known upper bound for both conjectures is 
\cite{Swa}
\[s(\M^n)\leq d(\M^n)\leq c2^n n^2\log n.\]
It is known that $s(\M^2)\leq 4$ for all planes \cite{SwaN}.
There are many two-dimensional spaces attaining $s(\M^2)=4$.
Some piece-wise $C^\infty$, elliptic examples were discovered in 
\cite{Alfetal}.
They are characterised in \cite{SwaN}.

We give a lower bound for $s(\M^n_Z)$ that is asymptotically correct up 
to a factor of $2$ by Theorem~\ref{thm1}.
\begin{theorem}\label{thm2}
For any $n\geq3$, $s(\M^n_Z)\geq\binom{n+1}{\lfloor (n+1)/2\rfloor}$.
Therefore, $s(\M^n_Z) = \Theta(2^n/\sqrt{n})$.
\end{theorem}
The proof is in Section~\ref{equilateral}.

The sharp upper bound for differentiable norms is 
$s(\M^n)\leq d(\M^n)\leq n+1$ \cite{LM, Swalp}.
It is well-known that $s(\ell_\infty^n)=d(\ell_\infty^n)= 2^n$ 
\cite{Morgan} and $s(\ell_1^n)=d(\ell_1^n)=2n$ 
(see also Section~\ref{equilateral}).
For the $\ell_p$ norm, $1<p<\infty$, we have 
$3\leq s(\ell_p^n)\leq d(\ell_p^n)\leq 7$ if $p>2, n\geq 2$, and 
$\min\{n,\frac{p}{(p-1)\ln 2}\}\leq s(\ell_p^n)\leq d(\ell_p^n)\leq 2^{p/(p-1)}$ 
if $1<p<2, n\geq 3$; see \cite{Swalp}, where more detailed estimates 
are obtained.

Conger \cite{Conger} showed that 
$s(\R^3, \norm{\cdot}_1+\lambda\norm{\cdot}_2)\geq 6$ for all 
$0<\lambda\leq 1$.
These norms are piece-wise $C^\infty$ and elliptic.
In \cite{Alfetal} it is shown that 
$s(\R^2, \norm{\cdot}_1+\lambda\norm{\cdot}_2)=4$ for all 
$0<\lambda\leq 2+\sqrt{2}$.
The value $\lambda=2+\sqrt{2}$ is sharp, since it follows from the 
results in \cite{SwaN} that 
$s(\R^2, \norm{\cdot}_1+\lambda\norm{\cdot}_2)=3$ for all 
$\lambda>2+\sqrt{2}$.
We generalise these results as follows.
\begin{theorem}\label{mcconj}
Let $\lambda>0$, and let 
$\M^n=(\R^n,\norm{\cdot}_1+\lambda\norm{\cdot}_2)$.
If $\lambda\leq 1$ then $s(\M^n)\leq d(\M^n)\leq 2n$.
If $\lambda\leq\sqrt{n}/(\sqrt{n}-1)$ then 
$d(\M^n)\geq s(\M^n)\geq 2n$.
\end{theorem}
The proof is in Section~\ref{elliptic}.
In this regard Conger made the following conjecture 
\cite[p.~128]{Morganbook}.
\begin{conjecture}[Conger]
For any piece-wise differentiable, elliptic $\M^n$ we have 
$s(\M^n)\leq2n$.
\end{conjecture}

Our results are based on a characterization of the vertex figures 
of nodes and of Steiner points in SMTs in arbitrary Minkowski spaces 
(Section~\ref{local}).
This characterization is found using the subdifferential calculus 
(Section~\ref{subdifferential}).
The characterization of vertex figures of nodes enables us to reduce 
the determination of $d(\M^n_Z)$ to a purely combinatorial problem 
in extremal finite set theory (Section~\ref{extremal}).
The characterization of vertex figures of Steiner points in the case 
of $\M^n_Z$ is not so easily reducible to combinatorics, hence the 
partial results of Theorem~\ref{thm2} (Section~\ref{equilateral}).
This situation can be compared to \cite{SwaN}, where the characterization 
of nodes in an arbitrary two-dimensional Minkowski space is much simpler 
than the characterization of Steiner points.
However, the norm $\norm{\cdot}_1+\lambda\norm{\cdot}_2$ is sufficiently 
simple so that the maximum degree of a Steiner point can be determined 
(Section~\ref{elliptic}).

In the next section we collect basic definitions used in the paper.

\section{Basic definitions}\label{definitions}
\subsection{Signed sets}
Let $[n]:=\{1,2,\dots,n\}$, and let $\CP[n]$ be the set of all subsets 
of $[n]$.
Let $\card{A}$ denote the number of elements in the finite set $A$.
We use the following standard notation for signed sets (see \cite{OM}).
A \emph{signed subset} $X$ of $[n]$ is a vector $X\in\{0,\pm 1\}^n$.
We denote the $i$th component of $X$ by $X(i)$ (reserving subscripts 
such as $X_i$ for indices).
The \emph{positive part} of $X$ is $X^+:=\{i:X(i)=+1\}$, the 
\emph{negative part} of $X$ is $X^-:=\{i:X(i)=-1\}$, the \emph{support} 
of $X$ is $\underline{X}:=X^+\cup X^-$, and the \emph{zero set} of $X$ 
is $X^0:=[n]\setminus\underline{X}$.
We also write $X=(X^+,X^-)$.
Thus any two disjoint subsets $A, B$ of $[n]$ determine a signed set 
$(A,B)$.
The natural partial ordering on signed sets $X\leq Y$ is defined by 
$X^+\subseteq Y^+$ and $X^-\subseteq Y^-$.
Two signed sets $X, Y$ are \emph{conformal} if 
$X^+\cap Y^-=\emptyset=X^-\cap Y^+$.
Denote the \emph{empty signed set} $(\emptyset,\emptyset)$ by $\emptyset$.

\subsection{Vector spaces}
We consider $\R^n$ to be the vector space of column vectors 
$\vx=(\vx(1),\dots,\vx(n))^\Tr$, with the standard basis $\ve_1,\dots,\ve_n$ 
defined by $\ve_i(j)=\delta_{ij}$.
The \emph{support} of a vector $\vx\in\R^n$ is the signed set 
$\supp(\vx):=(\supp^+(\vx),\supp^-(\vx))$, where 
$\supp^+(\vx)=\{i:\vx(i)>0\}$ and $\supp^-(\vx)=\{i:\vx(i)<0\}$.
Let $\vx^+, \vx^-$ be defined by $\vx^+(i)=\max\{\vx(i),0\}$ and 
$\vx^-(i)=\min\{\vx(i),0\}$.

Let $V\subseteq\R^n$ be a subspace, and $\pi:\R^n\to V$ a projection.
Define $\supp(V):=\{\supp(\vx):\vx\in V\}$.
This is the set of covectors of the oriented matroid of the vector 
configuration $\{\pi(\ve_i):i=1,\dots,n\}$ in $V$; cf.\ \cite[\S 2.2]{OM}.

The \emph{dual} $V^\ast$ of a vector space $V$ is the set of all linear 
functionals $\fhi:V\to\R$.
We identify $(\R^n)^\ast$ with the set of row vectors 
$\fhi=(\fhi(1),\dots,\fhi(n))$, with standard basis $\epsilon_i=\ve_i^\Tr$.
Let $V$ be a subspace of $\R^n$, and let $\pi:\R^n\to V$ be an orthogonal 
projection.
The quotient space $\R^n/\ker\pi$ can then be identified with $V$.
The dual of $\R^n/\ker\pi$ is in a canonical way isomorphic to the 
annihilator of $\ker\pi$, i.e., the subspace 
$\{\fhi:\ker\pi\subseteq\ker\fhi\}$ of $(\R^n)^\ast$.
We accordingly identify this space with $V^\ast$.
It is easily seen that with these identifications 
$\supp(V)=\supp(V^\ast)$ (as long as $\pi$ is orthogonal).

For any two vector spaces $V$ and $W$, $(V\times W)^\ast$ is canonically 
isomorphic to $V^\ast\times W^\ast$, since a functional 
$\fhi\in(V\times W)^\ast$ may be decomposed uniquely as 
$\fhi(\vx,\vy)=\psi(\vx)+\chi(\vy)$ for some $\psi\in V^\ast$ and 
$\chi\in V^\ast$; in fact, $\psi(\vx)=\fhi(\vx,\vo)$ and 
$\chi(\vy)=\fhi(\vo,\vy)$.
Conversely, given any $\psi\in V^\ast, \chi\in W^\ast$, their sum 
$\psi(\vx)+\chi(\vy)$ defines a functional on $V\times W$ which we 
denote $(\psi,\chi)$.

If $A$ and $B$ are subsets of a vector space $V$, we define their 
\emph{Minkowski sum} to be $A+B:=\{\va+\vb:\va\in A,\vb\in B\}$.

\subsection{Minkowski spaces}\label{mspaces}
A \emph{Minkowski space} is a finite-dimensional vector space $V$ 
together with a norm $\norm{\cdot}$ on $V$, denoted by 
$\M=(V,\norm{\cdot})$, or $\M^n$ to indicate that it is $n$-dimensional.
The \emph{unit ball} $B=B(\M)=\{\vx\in V:\norm{\vx}\leq 1\}$ determines 
the norm uniquely, and we also write $\M=(V,B)$.

The \emph{dual} of $\M=(V,\norm{\cdot})$ is 
$\M^\ast=(V^\ast,\norm{\cdot}^\ast)$, with the dual norm defined by
\[ \norm{\fhi}^\ast := \max\{\fhi(\vx):\norm{\vx}=1\}.\]
Denote the dual unit ball by $B^\ast$.
It is well-known that the dual of the dual is again the original space.

The $\ell_p$-norm on $\R^n$ is defined by
\[\norm{(\vx(1),\dots,\vx(n))}_p:=(\sum_{i=1}^n\abs{\vx(i)}^p)^{1/p}\]
if $p\geq 1$, and
\[\norm{(\vx(1),\dots,\vx(n))}_\infty:=\max\{\vx(i):i\in[n]\}.\]
The space $\ell_p^n$ is $(\R^n,\norm{\cdot}_p)=(\R^n,B_p^n)$.
The \emph{$n$-cube} is $B_\infty^n$, and the 
\emph{$n$-dimensional cross polytope} is $B_1^n$.
The spaces $\ell_\infty^n$ and $\ell_1^n$ are dual.

By the (Hahn-Banach) separation theorem, each non-zero $\vx\in \M^n$ 
has a \emph{norming functional}, i.e.\ a functional $\fhi$ such that 
$\norm{\fhi}^\ast=1$ and $\fhi(\vx)=\norm{\vx}$.
We denote the set of norming functionals of a non-zero $\vx$ by 
$\subdiff\vx$.
Each $\subdiff\vx$ is an \emph{exposed face} of $B^\ast$, i.e., the 
intersection of $B^\ast$ by some supporting hyperplane, in fact 
$\subdiff\vx=B^\ast\cap\{\fhi:\fhi(\vx)=\norm{\vx}\}$.
A norm is differentiable if and only if $\subdiff\vx$ is a singleton 
for each $\vx\neq\vo$ (and then the norm is $C^1$).

A norm is \emph{strictly convex} if for all linearly independent 
$\vx,\vy\in\M^n$ we have $\norm{\vx+\vy}<\norm{\vx}+\norm{\vy}$.
This is equivalent to there being no straight segment in the boundary 
of the unit ball (i.e., \emph{strict convexity} of the unit ball).
A norm is strictly convex if and only if the dual norm is differentiable.

A norm $\norm{\cdot}$ on $V$ is \emph{elliptic} if 
$\norm{\cdot}-\norm{\cdot}_2$ is still a norm for some Euclidean norm 
$\norm{\cdot}_2$ (i.e.\ for any identification of $V$ with $\R^n$).
We also say that the unit ball of an elliptic norm is \emph{elliptic}.
If the norm is $C^2$, this is equivalent to requiring that the boundary 
of the unit ball has positive inward curvature bounded away from $0$, 
where the curvature is defined by any Euclidean structure on the space.
A norm is elliptic if and only if its dual is $C^{1,1}$.
See \cite[\S3]{M4} for more on elliptic norms.

Note that in functional analysis a differentiable norm is also called 
a smooth norm, while in differential geometry the word smooth is 
usually reserved for $C^\infty$.
In differential geometry the unit ball of an elliptic norm is also 
called uniformly convex, which again has another meaning in Banach 
space geometry.
We therefore avoid the terms smooth and uniformly convex in this 
paper, and stick to the terminology of the previous paragraph.

If $V$ is a subspace of $\M=(W,\norm{\cdot})$, then the quotient 
space $\M/V$ is the space $W/V$ with norm 
$\norm{\vx+V}:=\inf\{\norm{\vx+\vv}:\vv\in V\}$.
If $V$ is a subspace of $\M=(\R^n,\norm{\cdot})$, a concrete presentation 
of $\M/V$ can be obtained by using the orthogonal projection 
$\pi:\R^n\to V^\perp$.
Then $\M/V$ is isometrically isomorphic to $(V^\perp,\pi(B))$, i.e., 
we project the unit ball $B$ orthogonally onto $V^\perp$ to obtain 
the norm.
The dual $(\M/V)^\ast$ is as before the annihilator of $\ker\pi=V$, 
and the dual norm is a restriction of the dual norm 
$\norm{\cdot}^\ast$ of $\M^\ast$, i.e., the dual is a subspace of 
$\M^\ast$, with the embedding given by the adjoint $\pi^\ast$ of the 
projection $\pi$.

\subsection{Zonotopes as unit balls}\label{zonotopes}
Denote the closed segment from $\va$ to $\vb$ by $[\va,\vb]$.
A \emph{zonotope} is a Minkowski sum of closed segments, equivalently, 
a zonotope is the projection of an $n$-cube onto a subspace.
The zonotope \emph{determined by} a subspace $V$ of $\R^n$  is the 
orthogonal projection $Z=\pi(B_\infty^n)$ of the $n$-cube onto $V$.
As explained in Section~\ref{mspaces}, $(V,Z)$ is isometric to the 
quotient $\ell_\infty^n/\ker\pi$,
and its dual $(V^\ast,Z^\ast)$ is isometric to the subspace 
$\im\pi^\ast$ of $\ell_1^n$.
As before, we identify $V^\ast$ with this subspace.
Then $Z^\ast=V^\ast\cap B_1^n$.
We need a description of the faces of $Z^\ast$.
Each proper face of $Z^\ast$ is the intersection of $V^\ast$ with a 
proper face of $B_1^n$, i.e., there is one for each signed set 
$X\in\supp(V^\ast)$, namely
\[ F'(X) := \{\fhi\in V^\ast:\norm{\fhi}_1=1, \supp(\fhi)\leq X\},\]
of dimension $\card{\underline{X}}-1$.
Since $X\mapsto F'(X)$ is order-preserving, we have that $F'(X)$ is a 
facet if and only if $X$ is a maximal covector in $\supp(V^\ast)$.
Two faces $F'(X)$ and $F'(Y)$ belong to a common facet if and only 
if $X$ and $Y$ are conformal.

We mention that the corresponding non-empty faces of $Z$ are indexed 
by the same set $\supp(V)=\supp(V^\ast)$.
Each such face equals some
\[ F(X):=\sum_{i\in X^+}\pi(\ve_i)-\sum_{j\in X^-}\pi(\ve_j)
+\sum_{k\in X^0}[-\pi(\ve_k),\pi(\ve_k)],\]
of dimension $\card{X^0}$ (see \cite[\S 2.2]{OM}).
The mapping $X\mapsto F(X)$ is order-reversing, and so $F(X)$ is a 
vertex if and only if $X$ is a maximal covector in $\supp(V^\ast)$, 
while two faces $F(X)$ and $F(Y)$ intersect if and only if $X$ and 
$Y$ are conformal.

We now consider the Minkowski space $\M^n_Z$ defined in the introduction.
Let $\pi$ be the orthogonal projection of $\R^{n+1}$ onto the hyperplane
\[ H := \{\vx\in\R^{n+1} : \sum_{i=1}^{n+1}\vx(i)=0\}.\]
This projection is along the vector $\vj:= (1,1,\dots,1)^\Tr$.
Thus $\M^n_Z$ is the quotient $\ell_\infty^{n+1}/\ker\pi$, which we 
identify with the $n$-dimensional space $H$ with unit ball 
$Z_n:=\pi(B_\infty^{n+1})$.
It is easily seen that the norm is given by
\[ \norm{\vx}_Z := \tfrac12(\max_{i\in[n+1]} \vx(i)
-\min_{i\in[n+1]} \vx(i)).\]
The dual space $\M_C^n:=(\M_Z^n)^\ast$ is the subspace
\[ H^\ast := \{\fhi\in(\R^{n+1})^\ast : 
\sum_{i=1}^{n+1}\fhi(\ve_i)=0\}\]
of $\ell_1^{n+1}$, with support
\[\CZ_n:=\supp(H^\ast)=\{X\in\{0,\pm 1\}^{n+1}:X^+\neq\emptyset, 
X^-\neq\emptyset\}\cup\{\emptyset\}.\]
Note that $F'(\emptyset)=\emptyset$, and all other $F'(X)$ are 
non-empty proper faces of $Z_n^\ast$.
The facets of $Z_n^\ast$ are all the $F'(X)$ with $X\in\{\pm1\}^{n+1}$, 
$X^\pm\neq\emptyset$.
In general,
\[ F'(X)=\{\fhi\in(\R^{n+1})^\ast:\sum_{i\in X^+}\fhi(\ve_i)=\tfrac12, 
\sum_{i\in X^-}\fhi(\ve_i)=-\tfrac12,\; \supp(\fhi)\leq X\}.\]
Also define $\overline{\CZ_n}:=\CZ_n\cup\{\one\}$, with $\one$ a symbol 
corresponding to the improper face $Z_n^\ast$.
We consider $\one\geq X$ for all $X\in\CZ_n$.
Thus $(\overline{\CZ_n},\leq)$ is the face lattice of $Z_n^\ast$.

It is easy to see that $Z_n^\ast$ is the difference body of the 
regular $n$-simplex \[\Delta_n:=\conv\{\epsilon_i:i\in[n+1]\}.\]
Thus $Z_2^\ast$ is a regular hexagon, and $Z_3^\ast$ a cuboctahedron.

To show that $Z_n$ is affinely equivalent to $\conv([0,1]^n\cup[-1,0]^n)$, 
let $\pi'$ be the (non-orthogonal) projection of $\R^{n+1}$ onto the 
hyperplane
\[ H' = \{\vx\in\R^{n+1}: \vx(n+1)=0\},\]
again along the vector $\vj$.
Then $\pi'|_{H}$ and $\pi|_{H'}$ are inverses of each other taking 
$Z_n$ to $\pi'(C^{n+1})$ and vice versa.
Identifying $H'$ with $\R^n$ in the obvious way, it is easily seen 
that $\pi'(C^{n+1})=\conv([0,2]^n\cup[-2,0]^n)$.

\section{Subdifferential calculus}\label{subdifferential}
In this section we summarise the basic results of the subdifferential 
calculus needed in this paper.
For more on subdifferential calculus, see \cite{Rockafellar} and 
\cite{Zalinescu}.
We only consider convex functions defined everywhere in 
finite-dimensional spaces, so many of the arguments simplify.
Although these results are all well-known in convex analysis, they are 
not easy to locate in the same reference, and the proofs given are 
usually based on a variety of other more general theorems.
For the convenience of the reader we give direct proofs of these lemmas.
We use the following version of the (weak) separation theorem 
\cite[Theorem~11.3]{Rockafellar}.
\begin{separationtheorem}
Let $C_1$ and $C_2$ be non-empty convex sets with disjoint relative 
interiors in a finite-dimensional real vector space.
Then there exists a hyperplane $H$ such that $C_1$ and $C_2$ are in 
opposite closed half spaces bounded by $H$.
\end{separationtheorem}

Let $V$ be any Minkowski space with dual $V^\ast$.
As usual, a function $f:V\to\R$ is \emph{convex} if 
$f(\lambda\vx+(1-\lambda)\vy)\leq\lambda f(\vx)+(1-\lambda)f(\vy)$ 
for all $\vx,\vy\in V$ and $0\leq\lambda\leq 1$.
Convex functions on $V$ are continuous on $V$.
A functional $\fhi\in V^\ast$ is a \emph{subgradient} of a convex 
function $f$ at the point $\va$ if for all $\vx\in V$,
\[ f(\vx)-f(\va) \geq\fhi(\vx-\va).\]
In particular, $o\in V^\ast$ is a subgradient of $f$ at $\va$ if 
and only if $f$ attains its minimum value at $\va$.
The \emph{subdifferential} of $f$ at $\va$, denoted by 
$\subdiff f(\va)$, is the set of all subgradients of $f$ at $\va$.
Note that for any $\fhi\in V^\ast$ and $c\in\R$, 
$\subdiff(\fhi+c)(\va)=\{\fhi\}$ for all $\va\in V$.
\begin{lemma}\label{l1}
The subdifferential of a convex function $f: V\to\R$ at any $\va\in V$ 
is non-empty, compact, and convex.
\end{lemma}
See \cite[Theorem~24.7]{Rockafellar}.
\begin{proof}
Since $\subdiff f(\va)=\bigcap_{\vx\in V}
\{\fhi\in V^\ast:\fhi(\vx-\va)\leq f(\vx)-f(\va)\}$,
which is an intersection of closed half spaces in $V^\ast$, we have 
that $\subdiff f(\va)$ is closed and convex.
If $\subdiff f(\va)$ is unbounded, a compactness argument shows that 
it contains a ray $\{\fhi_0+\lambda\fhi_1:\lambda\geq 0\}$.
Therefore, $f(\vx)-f(\va)\geq\fhi_0(\vx-\va)+\lambda\fhi_1(\vx-\va)$ 
for all $\vx\in V$ and $\lambda\geq 0$.
It follows that $\fhi_1(\vx-\va)\leq 0$ for all $\vx\in V$, a 
contradiction.

To show that $\subdiff f(\va)\neq\emptyset$ we use the separation 
theorem to obtain a hyperplane that separates the point 
$(\va,f(\va))$ and the \emph{epigraph} of $f$
\[ \{(\vx,y):\vx\in V, y\geq f(\vx)\}\subseteq V\times\R.\]
Thus there exists a non-zero functional 
$(\fhi, r)\in(V\times\R)^\ast=V^\ast\times\R$
such that $(\fhi, r)(\vx,y)\geq(\fhi, r)(\va,f(\va))$ for all $\vx\in V$ 
and $y\geq f(\vx)$.
Then $\fhi(\vx-\va)\geq r(f(\va)-y)$, which gives that $r\geq 0$, 
since $y$ can be arbitrarily large.
If $r=0$ this also gives $\fhi=o$, contradicting $(\fhi, r)\neq (o, 0)$.
Therefore, $r>0$, and $-\frac{1}{r}\fhi(\vx-\va)\leq f(\vx)-f(\va)$ 
for all $\vx\in V$, giving $-\frac{1}{r}\fhi\in\subdiff f(\va)$.
\end{proof}
\begin{lemma}[Moreau-Rockafellar Theorem]\label{l2}
Let $f_i:V\to\R$, $i=1,\dots,k$ be convex functions.
Then for all $\va\in V$,
\[ \subdiff(\sum_{i=1}^k f_i)(\va) = \sum_{i=1}^k\subdiff f_i(\va),\]
where the sum on the right is Minkowski addition.
\end{lemma}
See \cite[Theorem~23.8]{Rockafellar} and 
\cite[Theorem~2.8.7]{Zalinescu}.
\begin{proof}
The ``$\supseteq$''-inclusion is straightforward, but the 
``$\subseteq$''-inclusion needs the separation theorem 
(see the second proof on p.~224 of \cite{Rockafellar}).
We prove the case $k=2$, with the general case following by induction.

Let $f$ and $g$ be convex functions on $V$ and let 
$\fhi\in\subdiff(f+g)(\va)$.
Then
\begin{equation}\label{l2star}
f(\vx)-f(\va)+g(\vx)-g(\va)\geq\fhi(\vx-\va) \text{ for all }\vx\in V.
\end{equation}
The sets
\[ C_1=\{(\vx, r)\in V\times\R: f(\vx)-f(\va)-\fhi(\vx-\va)\leq r\}\]
and
\[ C_2=\{(\vx, r)\in V\times\R: g(\vx)-g(\va)\leq-r\}\]
are both convex and closed with non-empty interior, and by 
\eqref{l2star} their interiors are disjoint.
They have a common point $(\va, 0)$.
By the separation theorem there exists a non-zero 
$(\psi, s)\in(V\times\R)^\ast=V^\ast\times\R$ such that
$(\psi, s)(\vx,r)\geq(\psi, s)(\va,0)$ for all $(\vx,r)\in C_1$ and 
$(\psi, s)(\vx,r)\leq(\psi,s)(\va,0)$ for all $(\vx,r)\in C_2$.
This gives that for all $\vx\in V$,
\begin{equation}\label{l2starstar}
\psi(\vx-\va)\geq-sr\text{ for all } r\geq f(\vx)-f(\va)-\fhi(\vx-\va),
\end{equation}
and
\[ \psi(\vx-\va)\leq-sr\text{ for all } r\leq g(\va)-g(\vx).\]
Since $r$ can be arbitrarily large in \eqref{l2starstar}, $s\geq 0$.
If $s=0$, then $\psi(\vx-\va)=0$ for all $\vx\in V$, contradicting 
$(\psi, s)\neq (o, 0)$.
Therefore, $s>0$, and we obtain
\[ -\tfrac{1}{s}\psi(\vx-\va)\leq f(\vx)-f(\va)-\fhi(\vx-\va) \]
and
\[ \tfrac{1}{s}\psi(\vx-\va)\leq g(\vx)-g(\va).\]
This gives $\fhi-\frac{1}{s}\psi\in\subdiff f(\va)$ and 
$\frac{1}{s}\psi\in\subdiff g(\va)$.
Adding, we obtain $\fhi\in\subdiff f(\va)+\subdiff g(\va)$, as 
required.
\end{proof}

\begin{lemma}\label{l3}
Let $f_i:V_i\to\R$, $i=1,\dots,k$ be convex, and define 
$f:V_1\times\dots\times V_k\to\R$ by 
$f(\vx_1,\dots,\vx_k)=\sum_{i=1}^k f_i(\vx_i)$.
Then the subdifferential of the convex function $f$ is the 
Cartesian product
\[ \subdiff f(\va_1,\dots,\va_k) = \prod_{i=1}^k
\subdiff f_i(\va_i)\subseteq V_1^\ast\times\dots\times V_k^\ast.\]
\end{lemma}
See \cite[Corollary~2.4.5]{Zalinescu}.
\begin{proof}
By induction it is sufficient to consider the case $k=2$.
Let $f:V\to\R$ and $g:W\to\R$ be convex, define $F: V\times W\to\R$ 
by $F(\vx,\vy)=f(\vx)+g(\vy)$, and choose $\va\in V$, $\vb\in W$.

First choose $\fhi\in\subdiff f(\va)$ and $\psi\in\subdiff g(\vb)$.
Then $f(\vx)-f(\va)\geq\fhi(\vx-\va)$ for all $\va\in V$, and 
$g(\vy)-g(\vb)\geq\psi(\vy-\vb)$ for all $\vb\in W$.
Adding, we obtain 
$F(\vx,\vy)-F(\va,\vb)\geq (\fhi,\psi)(\vx,\vy)-(\fhi,\psi)(\va,\vb)$, 
which gives $(\fhi,\psi)\in\subdiff F(\va,\vb)$.

Now choose $(\fhi,\psi)\in\subdiff F(\va,\vb)$.
This means that for all $\vx\in V$ and $\vy\in W$,
\[ \fhi(\vx-\va)+\psi(\vy-\vb)\leq f(\vx)-f(\va)+g(\vy)-g(\vb).\]
Setting $\vy=\vb$ we obtain $\fhi\in\subdiff f(\va)$, and setting 
$\vx=\va$ gives $\psi\in\subdiff g(\vb)$.
\end{proof}

We have already defined the set of norming functionals $\subdiff\vx$ 
of a non-zero $\vx\in\M$.
We now extend this notation by defining $\subdiff\vo=B(\M^\ast)$.
The reason is the following:
\begin{lemma}\label{l4}
The subdifferential of the norm of $\M$ at $\va$ is given by 
$\subdiff\norm{\va}=\subdiff\va$.
\end{lemma}
See \cite[Corollary~2.4.16]{Zalinescu}.
\begin{proof}
By definition, $\fhi\in B(\M^\ast)$ if and only if 
$\fhi(\vx)\leq\norm{\vx}$ for all $\vx\in\M$, which is equivalent to 
$\fhi\in\subdiff\norm{\vo}$.
This takes care of the case $\va=\vo$.

If $\va\neq\vo$, then for any $\fhi\in\subdiff\va$ and any $\vx\in\M$,
\[ \fhi(\vx-\va)=\fhi(\vx)-\norm{\va}\leq\norm{\vx}-\norm{\va},\]
hence $\fhi\in\subdiff\norm{\va}$.
Conversely, suppose
\begin{equation}\label{l4star}
\norm{\vx}-\norm{\va}\geq\fhi(\vx-\va)\text{ for all }\vx\in\M.
\end{equation}
Setting $\vx=\vo$ in \eqref{l4star} we obtain $\fhi(\va)\geq\norm{\va}$, 
giving $\norm{\fhi}^\ast\geq 1$.
Setting $\vx=2\va$ in \eqref{l4star} we obtain $\norm{\va}\geq\fhi(\va)$, 
giving $\fhi(\va)=\norm{\va}$.
Now \eqref{l4star} gives
$\norm{\vx}\geq\fhi(\vx)$ for all $\vx\in\M$, implying 
$\norm{\fhi}^\ast=1$.
It follows that $\fhi$ is a norming functional of $\va$.
\end{proof}
The distance function $\rho(\vx,\vy):=\norm{\vx-\vy}$ is easily seen 
to be convex on $\M\times\M$ by the triangle inequality.
\begin{lemma}\label{l5}
For any $\va,\vb\in\M$,
\[ \subdiff\rho(\va,\vb)=\{(\fhi,-\fhi):\fhi\in\subdiff(\va-\vb)\}
\subseteq\M^\ast\times\M^\ast.\]
\end{lemma}
\begin{proof}
Let $\fhi\in\subdiff(\va-\vb)$.
By Lemma~\ref{l4} we have for any $\vx,\vy\in\M$,
\begin{align*}
 \norm{\vx-\vy}-\norm{\va-\vb}&\geq\fhi(\vx-\vy-\va+\vb)\\
&= \fhi(\vx-\va)-\fhi(\vy-\vb) = (\fhi,-\fhi)(\vx-\va,\vy-\vb),
\end{align*}
which gives $(\fhi,-\fhi)\in\subdiff\rho(\va,\vb)$.

Conversely, let $(\fhi,\psi)\in\subdiff\rho(\va,\vb)$.
Then, in particular, for any $\vx\in\M$,
\[\rho(\vx,\vx)-\rho(\va,\vb)\geq\fhi(\vx-\va)+\psi(\vx-\vb),\]
i.e., $(\fhi+\psi)(\vx)\leq\fhi(\va)+\psi(\vb)-\norm{\va-\vb}$.
Since the right-hand side is independent of $\vx$, we obtain 
$\fhi+\psi=o$.
\end{proof}

\section{The local Steiner problem in Minkowski spaces}\label{local}
Our characterization of vertex figures is in terms of so-called 
reduced Minkowski addition defined on subsets of the dual of a 
Minkowski space.
\subsection{Reduced Minkowski addition}\label{sect:reduced}
We define the \emph{reduced Minkowski sum} of two closed, convex 
subsets $C$ and $D$ of the unit ball of $(\M^n)^\ast$ to be
\[ C\rma D = \{\fhi+\psi: \fhi\in C, \psi\in D, 
\norm{\fhi+\psi}^\ast\leq 1\},\]
i.e.\ $C\rma D$ is the usual Minkowski sum $C+D$ intersected by 
the unit ball $B^\ast$ of $(\M^n)^\ast$.
Reduced Minkowski addition is clearly commutative but not associative, 
although it satisfies a weak associative law 
(see Proposition~\ref{weakassoc} below).
We now consider some elementary properties of this binary operation.
Consider a finite non-empty family $\Sigma=\{A_i:i\in I\}$ of 
operands, where $I\subseteq[k]$ and each $A_i$ is a closed 
convex subset of $B^\ast$.
We call $I$ the \emph{support} of $\Sigma$.
A \emph{parenthesization} of $\Sigma$ is defined to be a 
parenthesization, in the usual sense, of some ordering
\begin{equation}\label{parenthesization} A_{j(1)}\rma\dots\rma A_{j(i)}
\end{equation}
of $\Sigma$, where $j:[\card{I}]\to I$ is some bijection.
Stanley calls this a \emph{binary set bracketing} of $\Sigma$ 
\cite[p.~178]{Stanley}, and its enumeration is called 
Schr\"oder's third problem \cite{Schroder}.
We denote a parenthesization of $\Sigma$ by $\langle \Sigma\rangle$, 
and we also call the family $I$ the \emph{support} of 
$\langle \Sigma\rangle$.
We may define a parenthesization of $\Sigma$ recursively as follows.
\begin{itemize}
\item If $\Sigma$ is the singleton $\{A_i\}$, then 
$\langle \Sigma\rangle =A_i$ is the only parenthesization of $\Sigma$.
\item If $\Sigma=\Sigma_1\cup\Sigma_2$, with $\Sigma_1$ and 
$\Sigma_2$ disjoint, and $\langle \Sigma_i\rangle $ is a 
parenthesization of $\Sigma_i$ for each $i=1,2$, then both 
$\langle \Sigma\rangle =(\langle \Sigma_1\rangle 
\rma\langle \Sigma_2\rangle)$ and 
$\langle \Sigma\rangle =(\langle \Sigma_2\rangle 
\rma\langle \Sigma_1\rangle )$ are \emph{equivalent} 
parenthesizations of $\Sigma$.
\end{itemize}
More generally, two parenthesizations of $\Sigma$ are 
\emph{equivalent} if they can be transformed into each other using the 
commutative law on any subexpression.
This is an equivalence relation, and two equivalent parenthesizations 
clearly evaluate to the same set.
The collection of equivalence classes of parenthesizations of a set 
$\Sigma$ with support $I$ corresponds bijectively with the 
collection of abstract trees with $\card{I}+1$ leaves labelled by 
the elements of $\{0\}\cup I$, and with $\card{I}-1$ internal 
vertices of degree $3$.
We call such a tree a \emph{rooted abstract Steiner tree} on $I$.
The node $0$ is the root of the tree, node $i$ corresponds to $A_i$ 
for each $i\in I$, and each internal vertex corresponds to an 
instance of $\rma$ in the parenthesization corresponding to the tree.
Note that we do not distinguish between left and right branches 
(as is done in the well-known bijection between parenthesizations 
of ordered expressions and planted, trivalent plane trees, both 
counted by the Catalan numbers; see \cite[Exercise~6.19(b),(f)]{Stanley}).
Stanley calls this a \emph{binary total partition} 
\cite[Example~5.2.6]{Stanley} of the set $I$.
The explicit construction of the bijection for our variant is 
the following.
We denote the rooted abstract Steiner tree associated to 
$\langle \Sigma\rangle $ by $T_0\langle \Sigma\rangle $.
\begin{itemize}
\item If $\langle \Sigma\rangle =A_i$ then 
$T_0\langle \Sigma\rangle $ is the tree joining $0$ and $i$.
\item If $\langle \Sigma\rangle 
=(\langle \Sigma_1\rangle \rma\langle \Sigma_2\rangle)$ then 
$T_0\langle \Sigma\rangle $ is the tree obtained by identifying 
the vertices $0$ in $T_0\langle \Sigma_1\rangle $ and 
$T_0\langle \Sigma_2\rangle $ to a single vertex $x$ (corresponding 
to the instance of $\rma$ operating on $\langle \Sigma_1\rangle$ and 
$\langle \Sigma_2\rangle$), and joining $x$ to a new root $0$.
\end{itemize}
It is clear that two parenthesizations of $\Sigma$ are equivalent 
if and only if their associated rooted abstract Steiner trees are equal.

The solution to Schr\"oder's third problem, i.e., the number $a_k$ 
of equivalence classes of parenthesizations with support $[k]$, is 
the product of the first $k-1$ odd numbers:
\[a_k=\prod_{i=1}^{k-1}(2i-1).\]
This is seen as follows.
Since $a_k$ equals the number of rooted abstract Steiner trees on 
$[k]$, i.e., trees with $k+1$ leaves $\{0,1,\dots,k\}$ and $k-1$ 
internal vertices of degree $3$, we have $a_1=1$ and 
$a_{k+1}=(2k-1)a_k$, since we may subdivide any of the $2k-1$ edges 
of such a tree on $[k]$ and join the new vertex to $k+1$, to 
obtain such a tree on $[k+1]$.

We remark that the $k$th Catalan number can be derived from $a_k$.
Since there are $k-1$ $\rma$'s in a parenthesization of 
$\{A_i:i\in[k]\}$, there are $2^{k-1}$ parenthesizations in an 
equivalence class.
This gives $2^{k-1}a_k$ parenthesizations (in the usual sense) of 
the ordered expression $A_{\pi(1)}\rma\dots A_{\pi(k)}$, where 
$\pi$ is a permutation of $[k]$, taken over all permutations $\pi$.
Therefore, the number of parenthesizations of $A_1\rma\dots\rma A_k$ 
is $2^{k-1}a_k/k!$, which equals the $k$th Catalan number.
This derivation of the Catalan numbers is essentially the same as 
the classical combinatorial derivation of Rodrigues \cite{Rodrigues}.

We define an \emph{abstract Steiner tree} on $I$ to be a tree with 
set of leaves $I$, and with $\card{I}-2$ internal vertices of 
degree $3$.
The abstract Steiner tree associated to a rooted abstract Steiner 
tree on $I$ is obtained by contracting the root $0$ and its 
incident edge.
The abstract Steiner tree of a parenthesization 
$\langle \Sigma\rangle$ obtained in this way from 
$T_0\langle\Sigma\rangle$ is denoted by $T\langle \Sigma\rangle$.
The number of abstract Steiner trees on $[k]$ clearly equals 
$a_{k-1}$.

We call two parenthesizations \emph{weakly equivalent} if they can 
be transformed into each other using the commutative law on any 
subexpression or the associative law on the whole expression, i.e., 
a bracketed expression
\[((\langle \Sigma_1\rangle \rma\langle \Sigma_2\rangle )\rma
\langle \Sigma_3\rangle)\]
may be transformed into
\[(\langle \Sigma_1\rangle 
\rma(\langle \Sigma_2\rma\Sigma_3\rangle)),\]
and vice versa.
It is again clear that two parenthesizations $\Sigma$ and $\Sigma'$ 
are weakly equivalent if and only if their abstract Steiner trees 
are equal: $T\langle \Sigma\rangle =T\langle \Sigma'\rangle $.
It follows that the number of weak equivalence classes of 
parenthesizations on $[k]$ equals $a_{k-1}$.
The operation $\rma$ has the following weak associativity property.
\begin{proposition}\label{weakassoc}
Let $\langle \Sigma\rangle _1$ and $\langle \Sigma\rangle _2$ be two 
weakly equivalent parenthesizations of $\Sigma=\{A_i:i\in I\}$, where 
each $A_i$ is a closed convex subset of $B(\M^\ast)$.
Then $o\in\langle \Sigma\rangle _1$ if and only if 
$o\in\langle \Sigma\rangle _2$.
\end{proposition}
\begin{proof}
By the definition of weak equivalence, it is sufficient to show for any 
three $A_i\subseteq B(\M^\ast)$, $i=1,2,3$, that 
$o\in  A_1 \rma(A_2 \rma A_3)$ 
if and only if $o\in ( A_1 \rma A_2)
\rma A_3$.
However, since $A_i\subseteq B(\M^\ast)$, both statements are 
equivalent to $o\in A_1+A_2+A_3$.
\end{proof}

\subsection{The characterization}\label{characterization}
\begin{theorem}[Nodes]\label{nodes}
Let $N=\{\vp_0,\vp_1,\dots,\vp_k\}$ be a set of points in a Minkowski 
space $\M^n$.
Then the star joining $\vp_0$ to each $\vp_i$, $i\in[k]$, is an SMT of 
$N$ if and only if $\langle \Sigma\rangle \neq\emptyset$ for each 
parenthesization $\langle \Sigma\rangle $ of 
$\Sigma=\{\subdiff(\vp_i-\vp_0):i\in[k]\}$.
\end{theorem}
\begin{proof}
$\Rightarrow$:
Consider any parenthesization $\langle \Sigma\rangle $ of $\Sigma$ 
and its associated rooted abstract Steiner tree $T_0\langle \Sigma\rangle $.
We now turn this tree into a Steiner tree in $\M^n$.
Associate leaf $i$ with $\vp_i$ for each $i\in\{0\}\cup[k]$, and 
associate each of the $k-1$ internal vertices with a variable point 
$\vx_i\in\M^n$, $i\in[k-1]$.
Denote this Steiner tree by 
$T_0\langle \Sigma\rangle (\vx_1,\dots,\vx_{k-1})$.
Its length is
\[ L(\vx_1,\dots,\vx_{k-1}):=\ell(T_0\langle \Sigma\rangle 
(\vx_1,\dots,\vx_{k-1})).\]
The Steiner points $\vx_i$ in 
$T_0\langle \Sigma\rangle (\vx_1,\dots,\vx_{k-1})$ may coincide 
--- this results in the tree being in fact a contraction of 
$T_0\langle \Sigma\rangle $.
Note that $L:\M^n\times\dots\times\M^n\to\R$ is a convex function since
\[ L(\vx_1,\dots,\vx_{k-1})=\sum_{e\in E(T_0\langle \Sigma\rangle )}
\rho_e(\vx_1,\dots,\vx_{k-1}),\]
where $\rho_e(\vx_1,\dots,\vx_{k-1})=\norm{\va-\vb}$, with $\va$ and 
$\vb$ the two points in $N\cup\{\vx_1,\dots,\vx_{k-1}\}$ associated 
to the vertices incident to $e$.
%elements of $\{\vp_0,\dots,\vp_k,\vx_1,\dots,\vx_{k-1}\}$.
Each $\rho_e$ depends on only one or two of the variables $\vx_i$.
Since $T_0\langle \Sigma\rangle (\vp_0,\dots,\vp_0)$ is the star 
joining $\vp_0$ to all $\vp_i$, $i\in[k]$, which is an SMT by 
assumption, $L$ attains its minimum at $(\vp_0,\dots,\vp_0)$.
Thus
\begin{equation}\label{ast}
o\in\subdiff L(\vp_0,\dots,\vp_0) =\sum_{e\in E(T_0\langle \Sigma\rangle )}
\subdiff\rho_e(\vp_0,\dots,\vp_0)
\end{equation}
by Lemma~\ref{l2}.
By Lemmas~\ref{l3}, \ref{l4} and \ref{l5}, if $e=x_ix_j$, then 
\[ \subdiff\rho_e(\vp_0,\dots,\vp_0)
=\{(o,\dots,o,\underbrace{\fhi}_{\text{pos.\ } i},o,\dots,o,
\underbrace{-\fhi}_{\text{pos.\ } j},o,\dots,o) : \fhi\in B^\ast\}\]
while if $e=p_ix_j$, then 
\[ \subdiff\rho_e(\vp_0,\dots,\vp_0)=\{(o,\dots,o,
\underbrace{\fhi}_{\text{pos.\ } j},o,\dots,o) : 
\fhi\in \subdiff(\vp_i-\vp_0)\}.\]
By considering each coordinate $i\in[k-1]$ of \eqref{ast} we obtain a 
functional $\fhi_e\in(M^n)^\ast$ for each edge 
$e\in E(T_0\langle \Sigma\rangle )$ such that
\begin{itemize}
\item
$\fhi_e \in \begin{cases}
B^\ast & \text{ if $e=\vx_i\vx_j$ or $e=\vp_0\vx_i$,}\\
\subdiff(\vp_i-\vp_0) & \text{ if $e=\vp_i\vx_j$, $i\neq0$,}
\end{cases}$
\item and for each Steiner point $\vx_i$, $\fhi_e=\fhi_{f}+\fhi_{g}$, 
where $e$ is the incoming edge and $f, g$ the two outgoing edges of 
$\vx_i$, when the tree is directed away from the root $\vp_0$.
\end{itemize}
By induction on the definition of $T_0\langle \Sigma\rangle $ ($\equiv$ 
induction on subexpressions of $\langle \Sigma\rangle $) we obtain that 
$\fhi_e\in\langle \Sigma\rangle $, where $e=\vp_0\vx_i$ is the root edge.
This gives $\langle \Sigma\rangle \neq\emptyset$.

$\Leftarrow$: Consider any Steiner tree in $\M^n$ on 
$\{\vp_0,\dots,\vp_k\}$.
By subdividing points if necessary, we obtain a tree with leaves 
$\{\vp_0,\dots,\vp_k\}$ and with $k-1$ Steiner points $\vx_i$ of degree 
$3$, some of them possibly coinciding with each other or with the $\vp_i$.
This tree is the rooted abstract Steiner tree of some parenthesization 
$\langle \Sigma\rangle $ of $\Sigma$.
%It is given that this parenthesization is non-empty.
As in the ``$\Rightarrow$''-argument, we obtain that 
$\langle \Sigma\rangle \neq\emptyset$ implies that 
$o\in\subdiff L(\vp_0,\dots,\vp_0)$, i.e., $L$ attains its minimum at 
$(\vp_0,\dots,\vp_0)$, which is when the tree is 
$T_0\langle\Sigma\rangle(\vp_0,\dots,\vp_0)$, the star joining $\vp_0$ 
to the other $\vp_i$.
\end{proof}

\begin{theorem}[Steiner points]\label{steiner}
Let $N=\{\vp_1,\dots,\vp_k\}$ be a set of points in a Minkowski space 
$\M^n$.
Then the star joining $\vp_0$ to each $\vp_i$, $i\in[k]$, is an SMT of 
$N$ if and only if $o\in\langle \Sigma\rangle $ for each parenthesization 
$\langle \Sigma\rangle $ of $\Sigma=\{\subdiff(\vp_i-\vp_0):i\in[k]\}$, 
if and only if $o\in\langle \Sigma'\rangle + \subdiff(\vp_k-\vp_0)$ for 
each parenthesization $\langle \Sigma'\rangle $ of 
$\Sigma'=\{\subdiff(\vp_i-\vp_0):i\in[k-1]\}$.
\end{theorem}
\begin{proof}
Since any parenthesization $\langle \Sigma\rangle $ is weakly equivalent 
to $\langle \Sigma'\rangle \rma\subdiff(\vp_k-\vp_0)$ for some 
parenthesization of $\langle \Sigma'\rangle $, Proposition~\ref{weakassoc} 
gives the equivalence between the two conditions.
We next show that the first condition is necessary and sufficient.

$\Rightarrow$:
This is similar to the proof of Theorem~\ref{nodes}.
Consider any parenthesization $\langle \Sigma\rangle $ of $\Sigma$.
Turn its associated rooted abstract Steiner tree 
$T_0\langle \Sigma\rangle $ into an abstract Steiner tree 
$T\langle \Sigma\rangle $ by contracting the root.
We turn $T\langle \Sigma\rangle $ into a directed graph as follows.
Denote the (new) edge into which the root was contracted by $\overline{e}$, 
and give it both directions, denoting the two directed edges by 
$e^+$ and $e^-$.
Give all other edges of $T\langle \Sigma\rangle $ a single direction 
away from $\overline{e}$.
This tree becomes a Steiner tree in $\M^n$ as follows.
Associate leaf $i$ with $\vp_i$ for each $i\in[k]$, and associate 
each of the $k-2$ internal vertices with a variable point 
$\vx_i\in\M^n$, $i\in[k-2]$.
Denote this Steiner tree by 
$T\langle \Sigma\rangle (\vx_1,\dots,\vx_{k-2})$.
Its length is the convex function
\[ L(\vx_1,\dots,\vx_{k-2}):=\ell(T\langle \Sigma\rangle 
(\vx_1,\dots,\vx_{k-2})).\]
Again note that the Steiner points $\vx_i$ in 
$T\langle \Sigma\rangle (\vx_1,\dots,\vx_{k-2})$ may coincide, and then 
the tree is a contraction of $T\langle \Sigma\rangle $.
Again, $L$ attains its minimum at $(\vp_0,\dots,\vp_0)$.
Calculating the subdifferential coordinate-wise, we obtain a functional 
$\fhi_{\vec{e}}\in(M^n)^\ast$ for each directed edge $\vec{e}$ of 
$E(T\langle \Sigma\rangle )$ such that
\begin{itemize}
\item $\fhi_{\vec{e}} \in \begin{cases}
B^\ast & \text{ if $\vec{e}$ is incident with two Steiner points,}\\
\subdiff(\vp_i-\vp_0) 
& \text{ if $\vec{e}$ is incident with $p_i$, $i\neq0$,}
\end{cases}$
\item $\fhi_{e^+}=-\fhi_{e^-}$, and 
\item for each Steiner point $\vx_i$, 
$\fhi_{\vec{e}}=\fhi_{\vec{f}}+\fhi_{\vec{g}}$, where $\vec{e}$ is the 
incoming edge and $\vec{f}, \vec{g}$ the two outgoing edges of $\vx_i$, 
with the convention that we ignore the outgoing $e^{+}$ or $e^-$ if 
$\vx_i$ is incident with $\overline{e}$.
\end{itemize}
Write $\langle \Sigma\rangle =\langle \Sigma^+\rangle \rma\langle \Sigma^-\rangle $, 
where $e^\pm$ points to the subtree associated with 
$\langle \Sigma^\pm\rangle $.
Let $I^{\pm}$ be the support of $\langle \Sigma^{\pm}\rangle $, 
and for each $i\in[k]$ let $\fhi_i=\fhi_{\vec{e}}$, where 
$\vec{e}$ is incident with $\vx_i$.
Again by induction on subexpressions we obtain that 
$\fhi_{e^\pm}=\sum_{i\in I^{\pm}}\fhi_i\in\langle \Sigma^{\pm}\rangle $.
From $\fhi_{e^+}=-\fhi_{e^-}$ it follows that 
$o\in\langle \Sigma\rangle$.

$\Leftarrow$: Similar to the corresponding direction in the proof 
of Theorem~\ref{nodes}.
\end{proof}
By the discussion in Section~\ref{sect:reduced}, when applying 
Theorem~\ref{nodes} (Theorem~\ref{steiner}) there are $a_k$ 
($a_{k-1}$, respectively) parenthesizations to consider.

Note that in the  ``$\Leftarrow$''-directions of the above proofs we 
did not need the parts of Lemmas~\ref{l1} to \ref{l5} depending on 
the separation theorem, i.e., Lemma~\ref{l1} and the 
``$\subseteq$''-part of Lemma~\ref{l2}.
These directions of Theorems~\ref{nodes} and \ref{steiner} are used 
to obtain lower bounds for $d(\M^n)$ and $s(\M^n)$.
On the other hand, to obtain upper bounds we need the 
``$\Rightarrow$''-directions, where the separation theorem is needed.

Recall that if the norm is differentiable, then $\subdiff\norm{\vx}$ 
is a singleton whenever $\vx\neq\vo$.
This drastically simplifies the conditions in Theorems~\ref{nodes} 
and \ref{steiner} and we regain the following two results.

\begin{corollary}[Lawlor and Morgan \cite{LM}]
Let $N=\{\vp_1,\dots,\vp_k\}$ be a set of points in a Minkowski 
space $\M^n$ with differentiable norm.
Let $\fhi_i$ be the norming functional of $\vp_i-\vp_0$, $i\in[k]$.
Then the star joining $\vp_0$ to each $\vp_i$, $i\in[k]$, is an SMT of 
$N$ if and only if \[\sum_{i=1}^k\fhi=o\] and for each subset 
$I\subseteq[k]$, \[\norm{\sum_{i\in I}\fhi}^\ast\leq 1.\]
\end{corollary}

\begin{corollary}[\cite{Swalp}]
Let $N=\{\vp_0,\vp_1,\dots,\vp_k\}$ be a set of points in a Minkowski 
space $\M^n$ with differentiable norm.
Let $\fhi_i$ be the norming functional of $\vp_i-\vp_0$, $i=1,\dots,k$.
Then the star joining $\vp_0$ to each $\vp_i$, $i=1,\dots,k$, is an 
SMT of $N$ if and only if for each subset $I\subseteq[k]$, 
\[\norm{\sum_{i\in I}\fhi}^\ast\leq 1.\]
\end{corollary}

\section{Using extremal set theory to prove Theorem~\ref{thm1}}
\label{extremal}
In order to determine $d(\M_Z^n)$ we apply Theorem~\ref{nodes}.
It follows from the discussion in Section~\ref{zonotopes} that norming 
functionals are described as follows.
\begin{lemma}
For any non-zero $\vx\in\M_Z^n$, we have $\subdiff\vx=F'(X)$, 
where $X\in\CZ_n$ is the unique signed set such that $F(X)$ is the 
face of $Z_n$ which has $\frac{1}{\norm{\vx}_Z}\vx$ in its 
relative interior.
\end{lemma}
We thus have to determine for which families 
$\Sigma=\{X_i:i\in[k]\}\subseteq\CZ_n$ all parenthesizations 
$\langle F'(\Sigma)\rangle $ of $F'(\Sigma):=\{F'(X_i):i\in[k]\}$ 
are non-empty.
To this end, we define a commutative, non-associative binary operation 
$\rmul$ on $\overline{\CZ_n}$ as follows:
\begin{itemize}
\item $X\rmul\one=\one\rmul X=X$ for all $X\in\overline{\CZ_n}$, and 
\item for all $X,Y\neq\one$,
\[ X\rmul Y := \begin{cases}
(X^+,Y^-) & \text{ if $X^+\cap Y^-=\emptyset$ and $X^-\cap Y^+\neq\emptyset$,}\\
(Y^+,X^-) & \text{ if $X^+\cap Y^-\neq\emptyset$ and $X^-\cap Y^+=\emptyset$,}\\
\one      & \text{ if $X^+\cap Y^-\neq\emptyset$ and $X^-\cap Y^+\neq\emptyset$,}\\
\emptyset & \text{ if $X^+\cap Y^-=\emptyset$ and $X^-\cap Y^+=\emptyset$.}
\end{cases}\]
\end{itemize}
Extending $F'$ to $\overline{\CZ_n}$ by defining 
$F'(\one)=\{o\}\subset\M_C^n$, we obtain the following.
\begin{lemma}\label{morphism}
For all $X,Y\in\overline{\CZ_n}$,
\[ F'(X\rmul Y)\subseteq F'(X)\rma F'(Y).\]
Furthermore, if at least one of the conditions $X=\one$, $Y=\one$, 
$X^+\cap Y^-=\emptyset$ or $X^-\cap Y^+=\emptyset$ holds, then
\[ F'(X\rmul Y) = F'(X)\rma F'(Y).\]
\end{lemma}
\begin{proof}
We assume that $X,Y\neq\one$, otherwise the lemma is trivial.

Consider the first relation.
Since this is trivial if $X\rmul Y=\emptyset$, we assume that 
$X$ and $Y$ are not conformal, say that $X^+\cap Y^-\neq\emptyset$.
There are now two cases, depending on whether $X^-\cap Y^+$ is 
empty or not.
\begin{itemize}
\item
If $X^-\cap Y^+=\emptyset$, we have to show that 
$F'(Y^+,X^-)\subseteq F'(X)\rma F'(Y)$.
Choose any $i\in X^+\cap Y^-$.
Then for any $\fhi\in F'(Y^+,X^-)$ we have $\fhi(\ve_i)=0$.
By setting $\psi=\fhi^-+\frac12\epsilon_i\in F'(X)$ and 
$\chi=\fhi^+-\frac12\epsilon_i\in F'(Y)$, we obtain 
$\fhi=\psi+\chi\in F'(X)\rma F'(Y)$.
\item
If $X^-\cap Y^+\neq\emptyset$, we have to show that 
$F'(\one)\subseteq F'(X)\rma F'(Y)$.
Choose any $i\in X^+\cap Y^-$ and $j\in X^-\cap Y^+$.
Then $\frac12\epsilon_i-\frac12\epsilon_j\in F'(X)$ and 
$-\frac12\epsilon_i+\frac12\epsilon_j\in F'(Y)$, giving 
$o\in F'(X)\rma F'(Y)$.
\end{itemize}
This establishes the inclusion.

Now consider the equality.
Without loss assume $X^+\cap Y^-=\emptyset$.
Let $\fhi\in F'(X)$ and $\psi\in F'(Y)$ be such that 
$\norm{\fhi+\psi}_1\leq 1$.
Then for each $i\in X^+$ we have 
$0\leq\psi(\ve_i)\leq\fhi(\ve_i)+\psi(\ve_i)$ and for each $i\in Y^-$ 
we have $0\geq \fhi(\ve_i) \geq \fhi(\ve_i)+\psi(\ve_i)$.
In particular, $X^+\cap\supp^-(\fhi+\psi)=\emptyset$ and 
$Y^-\cap\supp^+(\fhi+\psi)=\emptyset$, and
\begin{align*}
1 & \geq \norm{\fhi+\psi}_1\\
&=\sum_{i\in\supp^+(\fhi+\psi)}(\fhi(\ve_i)+\psi(\ve_i))+
\sum_{i\in\supp^-(\fhi+\psi)}-(\fhi(\ve_i)+\psi(\ve_i))\\
& \geq \sum_{i\in X^+}(\fhi(\ve_i)+\psi(\ve_i))+
\sum_{i\in Y^-}-(\fhi(\ve_i)+\psi(\ve_i))\\
& = \frac12+\sum_{i\in X^+}\psi(\ve_i) + \frac12+
\sum_{i\in Y^-}-\fhi(\ve_i)\\
& \geq 1.
\end{align*}
Thus equality holds everywhere, giving $\norm{\fhi+\psi}_1=1$, 
$\supp^+(\fhi+\psi)\subseteq X^+$, and 
$\supp^-(\fhi+\psi)\subseteq Y^-$.
Therefore, $\fhi+\psi\in F'(X^+,Y^-)$.
If $X^-\cap Y^+\neq\emptyset$, then $(X^+,Y^-)=X\rmul Y$, and we 
are done.
Otherwise, $X^-\cap Y^+=\emptyset$, giving that $X$ and $Y$ are 
conformal.
Then $F'(X\rmul Y)=\emptyset$, and also $F'(X)\rma F'(Y)=\emptyset$, 
since $F'(X)$ and $F'(Y)$ belong to the same facet of $Z_n^\ast$.
\end{proof}

\begin{theorem}
Let $\Sigma=\{X_i:i\in[k]\}\subseteq\CZ_n\setminus\{\emptyset\}$.
Then the following are equivalent.
\begin{enumerate}\renewcommand{\theenumi}{\alph{enumi}}
\item Some parenthesization of $F'(\Sigma)$ is empty.\label{one}
\item Some parenthesization of $\Sigma$ \textup{(}with operation 
$\rmul$\textup{)} equals $\emptyset$.\label{two}
\item There exist indices $a,b,c,d\in[k]$ with 
$\{a,b\}\cap\{c,d\}=\emptyset$ such that 
$(X_a^+\cup X_c^+)\cap(X_b^-\cup X_d^-)=\emptyset$.\label{three}
\end{enumerate}
\end{theorem}
\begin{proof}
The implication \eqref{one}$\implies$\eqref{two} follows from 
Lemma~\ref{morphism}, since $F'(X)=\emptyset$ implies $X=\emptyset$.

\eqref{two}$\implies$\eqref{three}: Let 
$\langle \Sigma\rangle =\emptyset$.
By the definition of $\rmul$, for some subexpression 
$\langle \Sigma_1\rangle \rmul\langle \Sigma_2\rangle $ of 
$\langle \Sigma\rangle $ we have 
$\langle \Sigma_1\rangle \neq\emptyset\neq\langle \Sigma_2\rangle $ 
but $\langle \Sigma_1\rangle \rmul\langle \Sigma_2\rangle =\emptyset$.
If both $\Sigma_1$ and $\Sigma_2$ are singletons, say 
$\Sigma_1=\{X_i\}$ and $\Sigma_2=\{X_j\}$, then $X_i$ and $X_j$ 
are conformal, i.e., $X_i^+\cap X_j^-=\emptyset=X_j^+\cap X_i^-$, 
so we may take $a=b=i$ and $c=d=j$.

If $\Sigma_1$ is a singleton $\{X_a\}$ but $\Sigma_2$ is not, 
then by induction on the definition of $\rmul$ there exist 
$X_c,X_d\in\Sigma_2$ such that $\langle \Sigma_2\rangle =(X_c^+,X_d^-)$ 
is conformal with $X_a$, i.e., 
$X_a^+\cap X_d^-=\emptyset= X_c^+\cap X_a^-$.
Thus we may set $b=a$.
A similar argument takes care of the case where $\Sigma_2$ is a 
singleton.

If $\Sigma_1$ and $\Sigma_2$ are both not singletons, then again by 
induction on the definition of $\rmul$ there exist 
$X_a,X_b\in\Sigma_1$ and $X_c,X_d\in\Sigma_2$ such that 
$\langle \Sigma_1\rangle =(X_a^+,X_b^-)$ and 
$\langle \Sigma_2\rangle =(X_c^+,X_d^-)$ are conformal, giving 
\eqref{three}.

\eqref{three}$\implies$\eqref{one}: It is sufficient to find a 
parenthesization of a subfamily of $F'(\Sigma)$ that equals the 
empty set.
First consider the case where $a$, $b$, $c$, $d$ are distinct 
indices.
Without loss assume that 
$F'(X_a)\rma F'(X_b)\neq\emptyset\neq F'(X_c)\rma F'(X_d)$.
By Lemma~\ref{morphism} we have 
$F'(X_a)\rma F'(X_b)=F'(X_a\rmul X_b)$, hence 
$X_a\rmul X_b\neq\emptyset$, and $X_a\rmul X_b=(X_a^+,X_b^-)$.
Similarly, $F'(X_c)\rma F'(X_d)=F(X_c^+,X_d^-)$, and again by 
Lemma~\ref{morphism},
\begin{align*}
&\quad (F'(X_a)\rma F'(X_b))\rma(F'(X_c)\rma F'(X_d))\\
&= F'(X_a^+,X_b^-)\rma F'(X_c^+,X_d^-)\\
&= F'(\emptyset)=\emptyset.
\end{align*}
Similarly,
\begin{itemize}
\item if $a=b$ and $c=d$, then $F'(X_a)\rma F'(X_c)=\emptyset$,
\item if $a=b$ and $c\neq d$, then 
$F'(X_a)\rma(F'(X_c)\rma F'(X_d))=\emptyset$, and
\item if $a\neq b$ and $c=d$, then 
$(F'(X_a)\rma F'(X_b))\rma F'(X_c)=\emptyset$.\qedhere
\end{itemize}
\end{proof}

The above theorem together with Theorem~\ref{nodes} now gives the 
following.

\begin{corollary}
Let $P=\{\vp_i:i\in[k]\}$ be a family of points in
 $\M_Z^n\setminus\{\vo\}$, with $\norm{\vp_i}_Z^{-1}\vp_i$ 
in the relative interior of face $F(X_i)$ of $Z_n$.
Then the star connecting $P$ to $\vo$ is an SMT of 
$P\cup\{\vo\}$ if and only if there do not exist indices 
$a,b,c,d\in[k]$ with $\{a,b\}\cap\{c,d\}=\emptyset$ and 
$(X_a^+\cup X_c^+)\cap(X_b^-\cup X_d^-)=\emptyset$.

Consequently, the star connecting $P$ to $\vo$ is an SMT of 
$P\cup\{\vo\}$ if and only if for any 
$I\subseteq[k]$, $2\leq\card{I}\leq 3$, the star connecting 
$\{\vp_i:i\in I\}$ to $\vo$ is an SMT of 
$\{\vp_i:i\in I\}\cup\{\vo\}$.
\end{corollary}
This also implies Corollary~\ref{cor1}.
The problem of determining $d(\M_Z^n)$ has now been reduced to 
a problem in extremal finite set theory.
\begin{theorem}\label{thm11}
Let $\{X_i:i\in[k]\}$ be a family of signed sets from $\CZ_n$ 
such that there do not exist indices $a,b,c,d\in[k]$ with 
$\{a,b\}\cap\{c,d\}=\emptyset$ and 
$(X_a^+\cup X_c^+)\cap(X_b^-\cup X_d^-)=\emptyset$.
Then $k\leq\binom{n+2}{\lfloor(n+2)/2\rfloor}$ with equality 
if and only if all $X_i^0=\emptyset$, and either all 
$\card{X_i^+}\in\{\lfloor(n+1)/2\rfloor,\lfloor(n+1)/2\rfloor+1\}$, 
or all 
$\card{X_i^+}\in\{\lceil(n+1)/2\rceil-1,\lceil(n+1)/2\rceil\}$.
\end{theorem}
\begin{proof}
It is easily seen that the hypothesis is equivalent to the 
following statement:
For all families of sets 
$\{Y_i:i\in[k]\}\subseteq\CP[n+1]\setminus\{\emptyset,[n+1]\}$ 
such that $X_i^+\subseteq Y_i\subseteq[n+1]\setminus X_i^-$,
\begin{equation}\label{star}
\begin{aligned}
&\text{there do not exist indices $a,b,c,d\in[k]$}\\
&\text{with $\{a,b\}\cap\{c,d\}=\emptyset$ and 
$Y_a\cup Y_c\subseteq Y_b\cap Y_d$.}
\end{aligned}
\end{equation}
Property \eqref{star} is in turn equivalent to the following 
three conditions:
\begin{align}
&\text{all $Y_i$ are distinct sets,}\label{distinct}\\
&\text{there do not exist distinct $a,b,c\in[k]$ with 
$Y_a\subseteq Y_b\subseteq Y_c$, and}\label{3chain}\\
&\text{there do not exist distinct $a,b,c,d\in[k]$ with 
$Y_a\cup Y_b\subseteq Y_c\cap Y_d$.}\label{butterfly}
\end{align}
By a well-known generalisation of Sperner's theorem due to 
Erd\H{o}s \cite{Erdos, Anderson}, conditions \eqref{distinct} 
and \eqref{3chain} on their own already give the sharp upper 
bound $k\leq \binom{n+2}{\lfloor(n+2)/2\rfloor}$ with equality 
exactly when $\{Y_i:i\in[k]\}$ consists of two of the largest 
levels in $\CP[n+1]$.
Because of this rigidity, it easily follows that in the case 
of equality all $X_i^0=\emptyset$, finishing the proof.
\end{proof}

Theorem~\ref{thm1} is now proved.
We remark that conditions~\eqref{distinct} and 
\eqref{butterfly} on their own give the same upper bound, 
but with an extra case of equality when $n+1=4$; see \cite{BKS}.

\section{A geometric formulation of Sperner's theorem}
\label{equilateral}
The next lemma follows from an observation of Moore \cite{GP}.
\begin{lemma}\label{moore}
If the unit ball of $\M^n$ contains $k$ points on its boundary 
such that the distance between any two equals $2$, then 
$s(\M^n)\geq k$.
\end{lemma}
From this lemma together with Lemma~\ref{polyhedral} it 
immediately follows that $s(\ell_\infty^n)=d(\ell_\infty^n)=2^n$ 
and $s(\ell_1^n)=d(\ell_1^n)=2n$.
We apply this lemma to $\M_Z^n$.
Recall that a vertex of $Z_n$ equals some $F(X)$, where 
$X^0=\emptyset$ and $X^\pm\neq\emptyset$.
Thus $X$ can be identified with the set $X^+\subseteq[n+1]$, 
not equal to $\emptyset$ or $[n+1]$.
We then have for any distinct 
$X^+,Y^+\in\CP[n+1]\setminus\{\emptyset,[n+1]\}$ that
\[\norm{F(X^+)-F(Y^+)}_Z=
\begin{cases}
1 & \text{if $X^+\subset Y^+$ or $Y^+\subset X^+$,}\\ 
2 & \text{if $X^+\not\subseteq Y^+$ and $Y^+\not\subseteq X^+$.}
\end{cases}\]
It follows that a set $\{F(X_i):i\in I\}$ of vertices of $Z_n$ 
are all at pairwise distance $2$ if and only if $\{X_i^+:i\in I\}$ 
is an antichain.
It follows from Lemma~\ref{moore} that 
$s(\M^n_Z)\geq \binom{n+1}{\lfloor(n+1)/2\rfloor}$.
This establishes Theorem~\ref{thm2}.

By Sperner's theorem \cite{Anderson} we cannot do better: 
The largest set of vertices at pairwise distance $2$ has size 
$\binom{n+1}{\lfloor(n+1)/2\rfloor}$.
We still cannot do better even if we consider arbitrary points 
on the boundary of $Z_n$.
For any boundary point $\vp$ of $Z_n$ its support 
$X_{\vp}:=\supp(\vp)$ describes the unique face $F(X_{\vp})$ 
which contains $\vp$ in its relative interior.
Two boundary points $\vp$ and $\vq$ are at distance $2$ if and only 
if there exist parallel supporting hyperplanes at $\vp$ and 
$\vq$ with $Z_n$ in between.
This in turn is equivalent to $\vp$ and $-\vq$ being contained in 
the same facet $F(X)$ of $Z_n$.
Note that $F(X)$ is a facet of $Z_n$ if and only if $\card{X^\pm}=1$.
Therefore, $\norm{\vp-\vq}_Z=2$ if and only if 
$X_{\vp}^+\cap X_{\vq}^-\neq\emptyset$ and 
$X_{\vp}^-\cap X_{\vq}^+\neq\emptyset$.
Now let $P=\{\vp_i:i\in I\}$ be a set of boundary points of $Z_n$ 
at pairwise distance $2$, let $X_i=\supp(\vp_i)$, and choose 
$Y_i\in\CP[n+1]\setminus\{\emptyset,[n+1]\}$ such that 
$X_i^+\subseteq Y_i\subseteq[n+1]\setminus X_i^-$.
Then it follows that the $Y_i$ are all distinct sets, and form 
an antichain.
As before we have by Sperner's theorem that 
$\card{P}\leq\binom{n+1}{\lfloor(n+1)/2\rfloor}$, with equality if 
and only if all $\vp_i$ are vertices and all $X_i^+$ have the same 
cardinality, either $\lfloor(n+1)/2\rfloor$ or $\lceil(n+1)/2\rceil$.
We have shown the following.
\begin{proposition}
The largest number of unit vectors in $\M_Z^n$ at pairwise distance 
$2$ is $\binom{n+1}{\lfloor(n+1)/2\rfloor}$. \qed
\end{proposition}
An even more general result would be the following.
\begin{conjecture}
The largest size of an equilateral set in $\M_Z^n$ equals 
$\binom{n+1}{\lfloor(n+1)/2\rfloor}$.
\end{conjecture}
This conjecture is known to hold for $n=3$ \cite{SS} 
(and is easy for $n=2$).

\section{Perturbations of $\ell_1^n$}\label{elliptic}
The following lemma is standard.
\begin{lemma}\label{ucdual}
Let $\norm{\cdot}_a$ and $\norm{\cdot}_b$ be any norms on $\R^n$.
Then the dual of $\M^n=(\R^n,\norm{\cdot}_a+\norm{\cdot}_b)$ 
is isometrically isomorphic to $((\R^n)^\ast, B_a^\ast+B_b^\ast)$, 
where $B_a^\ast$ and $B_b^\ast$ are the unit balls of the dual 
norms $\norm{\cdot}_a^\ast$ and $\norm{\cdot}_b^\ast$.
\end{lemma}

\begin{proof}[Proof of Theorem~\ref{mcconj}]
Let $\norm{\cdot}=\norm{\cdot}_1+\lambda\norm{\cdot}_2$.
By Lemma~\ref{ucdual} we may identify the dual unit ball 
$B^\ast$ with 
$\sum_{i=1}^n[-\epsilon_i,\epsilon_i]+\lambda B_2^\ast$.
Recall that each $\subdiff\norm{\vx}$, $\vx\neq\vo$, is an exposed 
face of $B^\ast$.
Any exposed face of $B^\ast$ equals a proper non-empty face $F$ of 
the cube $\sum_{i=1}^n[-\epsilon_i,\epsilon_i]$ translated by a 
functional $\fhi$ with $\norm{\fhi}_2^\ast=\lambda$, such that 
$F$ and $\fhi$ have the same sign in the following sense:
if $F=\sum_{i\in X^+}\epsilon_i-\sum_{i\in X^-}\epsilon_i%
+\sum_{i\in X^0}[-\epsilon_i,\epsilon_i]$ for some signed set 
$X\neq\emptyset$, then $X=\supp(\fhi)$.

We first show $d(\R^n,\norm{\cdot})\leq 2n$ if $\lambda\leq 1$.
Let $\vp_1,\dots,\vp_k\in\M^n$, $\vp_i\neq\vo$, and suppose that 
the star joining $\vo$ to all $\vp_i$ is an SMT of 
$\{\vo\}\cup\{\vp_i:i\in[k]\}$.
Then by Theorem~\ref{nodes} all restricted Minkowski sums of the 
$\subdiff\norm{\vp_i}$ must be non-empty.
Let $\subdiff\norm{\vp_i}=F_i+\fhi_i$ as above, with $X_i$ the 
corresponding signed set.
Suppose $k>2n$.
Since all $X_i\neq\emptyset$, it follows from the pigeon-hole 
principle that for some two indices $i,j$, $X_i$ and $X_j$ have a 
common element of the same sign, say $1\in X_1^+\cap X_2^+$.
Then, since $F_1$ and $F_2$ are both contained in the hyperplane 
$\{\chi:\chi(\ve_1)=1\}$, it follows that 
$\subdiff\norm{\vp_1}+\subdiff\norm{\vp_2}$ is contained in the 
hyperplane $\{\chi:\chi(\ve_1)=2+\fhi_1(\ve_1)+\fhi_2(\ve_1)\}$.
However, $B^\ast$ is contained in the slab bounded by 
$\pm\{\chi:\chi(\ve_1)=1+\lambda\}$.
Since $2+\fhi_1(\ve_1)+\fhi_2(\ve_1)>2\geq1+\lambda$, we obtain 
$\subdiff\norm{\vp_1}\rma\subdiff\norm{\vp_2}=\emptyset$, a 
contradiction.
Therefore, $k\leq 2n$.

\smallskip
We now prove that $s(\R^n,\norm{\cdot})\geq 2n$ for 
$\lambda\leq \sqrt{n}/(\sqrt{n}-1)$ by showing that the star joining 
$\vo$ to all $\pm\ve_i$, $i\in[n]$, is an SMT of 
$\{\pm\ve_i:i\in[n]\}$.
This is trivial for $n=1$, so we assume from now on that $n\geq 2$.
We have
\[ \subdiff\norm{\ve_i} = E_i := \sum_{\substack{j=1\\j\neq i}}^n
[-\epsilon_j,\epsilon_j]+(1+\lambda)\epsilon_i.\]
By Theorem~\ref{steiner} it is sufficient to prove that for each 
parenthesization $\langle\Sigma\rangle$ of
\[\Sigma:=\{\pm E_i:i\in[n-1]\}\cup\{E_n\}\] we have 
$o\in\langle\Sigma\rangle-E_n$, or equivalently, 
$E_n\cap\langle\Sigma\rangle\neq\emptyset$.

Write
\[\langle\Sigma\rangle=\langle\Sigma_1\rangle\rma\langle\Sigma_2\rangle
=(\langle\Sigma_1\rangle+\langle\Sigma_2\rangle)\cap B^\ast\]
where $\Sigma_1\cup\Sigma_2=\Sigma$, $\Sigma_1\cap\Sigma_2=\emptyset$, 
and $E_n\in\Sigma_1$.
Since $E_n\subseteq B^\ast$, we only have to prove that 
$E_n\cap(\langle\Sigma_1\rangle+\langle\Sigma_2\rangle)\neq\emptyset$.
We now replace the operation $\rma$ in $\langle\Sigma_i\rangle$ 
by $\rma'$, where
\[C\rma'D:=(C+D)\cap[-1-\lambda/\sqrt{n},1+\lambda/\sqrt{n}]^n.\]
Then $C\rma'D\subseteq C\rma D$, since 
$[-1-\lambda/\sqrt{n},1+\lambda/\sqrt{n}]^n\subseteq B^\ast$.
Denoting the parenthesizations with respect to $\rma'$ by 
$\langle\Sigma_i\rangle'$, it is sufficient to prove that 
$E_n\cap(\langle\Sigma_1\rangle'+\langle\Sigma_2\rangle')\neq\emptyset$.
Since all $E_i$ as well as $[-1-\lambda/\sqrt{n},1+\lambda/\sqrt{n}]^n$ 
are Cartesian products, we may show this coordinatewise.
For $C,D\subseteq\R$, let 
\[C\rma'' D:=(C+D)\cap[-1-\lambda/\sqrt{n},1+\lambda/\sqrt{n}],\] 
and for any family $\Sigma$ of subsets of $\R^n$, let $\pi_i(\Sigma)$ 
denote the family $\{\pi_i(A):A\in\Sigma\}$ of subsets of $\R$, 
where $\pi_i:\R^n\to\R$ is the $i$th coordinate projection.
We have to show the following:
\begin{equation}\label{first}
1+\lambda \in\langle\pi_n(\Sigma_1)\rangle''
+\langle\pi_n(\Sigma_2)\rangle'',
\end{equation}
and for all $i\in[n-1]$,
\begin{equation}\label{second}
[-1,1]\cap\langle\pi_i(\Sigma)\rangle''\neq\emptyset,
\end{equation}
where the parenthesizations are with respect to $\rma''$.
First note the following.
\begin{claim}\label{claimone}
Any parenthesization with respect to $\rma''$ of one or more sets 
all equal to $[-1,1]$ contains $[-1,1]$.
\end{claim}
It follows from Claim~\ref{claimone} and induction that
\begin{claim}\label{claimtwo}
Any parenthesization with respect to $\rma''$ of two or more sets 
all but one equal to $[-1,1]$, and the remaining set equal to 
$\{1+\lambda\}$, contains $[\lambda,1+\lambda/\sqrt{n}]$.
\end{claim}
The interval $[\lambda,1+\lambda/\sqrt{n}]$ is non-empty since 
$\lambda\leq\sqrt{n}/(\sqrt{n}-1)$ by hypothesis.

Since $\pi_n(\Sigma_1)$ consists of $[-1,1]$'s (perhaps none) 
and a $\{1+\lambda\}$, we obtain from Claim~\ref{claimtwo} that 
either 
$\langle\pi_n(\Sigma_1)\rangle''\supseteq[\lambda,1+\lambda/\sqrt{n}]$ 
or $\langle\pi_n(\Sigma_1)\rangle''=\{1+\lambda\}$.
Similarly, since $\pi_n(\Sigma_2)$ consists only of $[-1,1]$'s, 
we have $\langle\pi_n(\Sigma_2)\rangle''\supseteq[-1,1]$ by 
Claim~\ref{claimone}, and \eqref{first} follows.

Let $i\in[n-1]$.
Then $\pi_i(\Sigma)$ consists of a $\{1+\lambda\}$, a $\{-1-\lambda\}$, 
and $2n-3$ $[-1,1]$'s.
The parenthesization $\langle\pi_i(\Sigma)\rangle''$ has a unique 
subexpression 
$\langle\pi_i(\Sigma_0)\rangle''=\langle\pi_i(\Sigma_1)\rangle''%
\rma''\langle\pi_i(\Sigma_2)\rangle''$ such that 
$\{1+\lambda\}\in\pi_i(\Sigma_1)$ and 
$\{-1-\lambda\}\in\pi_i(\Sigma_2)$.
By Claim~\ref{claimtwo}, either 
$\langle\pi_i(\Sigma_1)\rangle''\supseteq[\lambda,1+\lambda/\sqrt{n}]$ 
or $\langle\pi_i(\Sigma_1)\rangle''=\{1+\lambda\}$.
Similarly, either 
$\langle\pi_i(\Sigma_2)\rangle''\supseteq[-1-\lambda/\sqrt{n},-\lambda]$ 
or $\langle\pi_i(\Sigma_2)\rangle''=\{-1-\lambda\}$.
This gives four cases, namely
\begin{itemize}
\item $\langle\pi_i(\Sigma_0)\rangle''%
\supseteq[-1+\lambda-\lambda/\sqrt{n},1-\lambda+\lambda/\sqrt{n}]$,
\item $\langle\pi_i(\Sigma_0)\rangle''%
\supseteq[-1,-\lambda+\lambda/\sqrt{n}]$,
\item $\langle\pi_i(\Sigma_0)\rangle''%
\supseteq[\lambda-\lambda/\sqrt{n},1]$,
\item $\langle\pi_i(\Sigma_0)\rangle''\supseteq\{0\}$.
\end{itemize}
In all four cases it follows that 
$\langle\pi_i(\Sigma_0)\rangle''\cap[-1,1]\neq\emptyset$ 
(the first three cases because $\lambda\leq1/(1-/\sqrt{n})$), 
and we obtain \eqref{second} from the following claim, which is 
proved similarly to Claim~\ref{claimtwo}.
\begin{claim}\label{claimthree}
Any parenthesization with respect to $\rma''$ of two or more sets 
all but one equal to $[-1,1]$, and the remaining set having 
non-empty intersection with $[-1,1]$, has non-empty intersection 
with $[-1,1]$.\qedhere
\end{claim}
\end{proof}

\section*{Acknowledgements}
I would like to thank Nic van Rensburg for drawing my attention to 
the subdifferential calculus, as well as Frank Morgan and 
Mark Conger for their remarks and encouragement.

\end{document}